\documentclass[reqno, 12pt]{amsart}
\usepackage{amsfonts}
\usepackage{amssymb,amsmath}
\usepackage{amsthm}
\usepackage{upref}
\usepackage{enumerate}
\usepackage{bbm}

\usepackage{amscd}

\makeatletter
\def\LaTeX{\leavevmode L\raise.42ex
    \hbox{\kern-.3em\size{\sf@size}{0pt}\selectfont A}\kern-.15em\TeX}
 \makeatother

 \newcommand{\1}{\mathbbm{1}}

\numberwithin{equation}{section}

\newtheorem{lemma}{Lemma}[section]
\newtheorem{theorem}[lemma]{Theorem} 
\newtheorem{corollary}[lemma]{Corollary}
\newtheorem{proposition}[lemma]{Proposition}

\theoremstyle{definition}

\newtheorem{example}[lemma]{Example}

\newtheorem{remark}[lemma]{Remark}

  \newcommand{\e}{\eqref}

\newcommand{\q}{\quad}

\newcommand{\wt}{\widetilde}

\newcommand{\la}{\langle}
\newcommand{\ra}{\rangle}
\newcommand{\z}{\zeta}
\newcommand{\ov}{\overline}
 \renewcommand{\d}{\delta}
 
 \newcommand{\cl}{\operatorname{clos}}
 
   \newcommand{\sgn}{\operatorname{sgn}}

\renewcommand\Re{\operatorname{Re}}

\newenvironment{pf}{\begin{proof}}{\end{proof}}

\def\qqq{\mathrel{\subset\mkern-15mu\lower.38ex\hbox{${\scriptscriptstyle\rightarrow}$}}}

\let\cal\mathcal

\let\Bbb\mathbb

\begin{document}

\title [Hankel  and Toeplitz operators] {   Hankel  and Toeplitz operators:    continuous and discrete representations  }
 
\author{ D. R. Yafaev}
\address{ IRMAR, Universit\'{e} de Rennes I\\ Campus de
  Beaulieu, 35042 Rennes Cedex, FRANCE}
\email{yafaev@univ-rennes1.fr}


\begin{abstract}
We   find a relation guaranteeing that Hankel operators   realized
  in the space of sequences $\ell^2 ({\Bbb Z}_{+}) $    and in the space    of functions $L^2 ({\Bbb R}_{+}) $ are unitarily equivalent. This allows us to obtain exhaustive spectral results for two classes of unbounded Hankel operators    in the space $\ell^2 ({\Bbb Z}_{+}) $ generalizing in different directions the classical Hilbert matrix.  We also discuss a link between   representations of Toeplitz operators in the spaces $\ell^2 ({\Bbb Z}_{+}) $ and $L^2 ({\Bbb R}_{+}) $.
       \end{abstract}

\maketitle


\section{Introduction}  

{\bf 1.1.}
This paper is based on the talk given by  the  author at the conference "Spectral Theory and Applications"   held in May 2015 in Krakow. So, it is somewhat eclectic. 
Our aim is to discuss various properties Hankel  and Toeplitz (known also as Wiener-Hopf) operators.
 We refer to the books \cite{Bo, GGK, NK, Pe, RoRo} for basic information on  
these classes of operators.

Our main goal is to describe a relation between discrete and continuous  representations of   Hankel  and Toeplitz operators in a sufficiently consistent way and to draw spectral consequences from this relation. We do not suppose that  operators are bounded, and so  we are naturally led   to work with quadratic forms and distributional integral kernels.   As is well known,  the discrete  (in the space $\ell^2 ({\Bbb Z}_{+}) $) and continuous (in the space $L^2 ({\Bbb R}_{+}) $)    representations are linked by the Laguerre transform. For bounded operators, this yields the unitary equivalence
of the corresponding operators in $\ell^2 ({\Bbb Z}_{+}) $ and $L^2 ({\Bbb R}_{+}) $. However in singular cases their equivalence may be lost because the natural domains of the quadratic forms in discrete and continuous  representations are not linked by the Laguerre transform.  As show simple examples, the continuous representation seems to be more    general. Passing to the Fourier transforms, one can also realize   discrete   Hankel  and Toeplitz operators in the Hardy space ${\Bbb H}_{+}^2 ({\Bbb T})$ of functions analytic in the unit circle  ${\Bbb T}$ and continuous   operators in the Hardy space ${\Bbb H}_{+}^2 ({\Bbb R})$ of functions analytic in the upper half-plane;  see, e.g., the book \cite{Hof}, for the precise  definition of these spaces.

Section~2 is of a preliminary nature. We first consider the discrete $A$ and the continuous $\bf A$  convolution  operators  in the spaces  $\ell^2  ({\Bbb Z} )$  and $L^2 ({\Bbb R} ;dx)$, respectively. Of course the Fourier transform allows one to reduce these operators to the multiplications $B$ and   $\bf B$ in the spaces $L^2 ({\Bbb T})$ and    
$L^2 ({\Bbb R}; d\lambda)$. The  operators  $B$ and   $\bf B$ are obviously related by a change of variables. This yields a link between the operators  $A$ and   $\bf A$ which is given by  the Laguerre transform. Our main goal in this section is to discuss explicit formulas relating matrix elements of $A$ and the integral kernel of $\bf A$.
Then, we apply these results to Toeplitz operators realized in the spaces $\ell^2  ({\Bbb Z}_{+} )$ and $L^2 ({\Bbb R}_{+} )$.  We are aiming at   a systematic presentation of known results but insist upon the case of unbounded operators. Moreover,  some formulas, for example,   \e{eq:K3b} and \e{eq:Kmod}, are perhaps new.

\medskip

{\bf 1.2.}
In Section~3, we pass  to the main subject of this paper, to  Hankel    operators.
  We recall that   Hankel  operators $H$    are defined in the space  $\ell^2  ({\Bbb Z_{+}})$ by the formula
  \begin{equation}
(H f )_{n}=\sum_{m\in {\Bbb Z}_{+}} a_{n+m} f_{m}, \q  f=\{ f_{n}\}_{n\in{\Bbb Z}_{+}}.
\label{eq:B1h}\end{equation}
Similarly,  Hankel operators $ {\bf H}  $ act in the space $L^2 ({\Bbb R}_{+})$ by the formula 
    \begin{equation}
({\bf H} {\bf f} )(t) = \int_{{\Bbb R}_{+}} {\bf a}(t+s) {\bf f}(s)  ds .
\label{eq:B6h}\end{equation}
Note that spectral properties of the operators $H$ are determined by the behavior of their matrix elements $a_{n} $ as $n\to\infty$ while, as far as the operators $\bf H$ are concerned, both the behavior of integral  kernels ${\bf a}(t)$ as $t\to \infty$ and $t\to 0$ as well as their local singularities at finite points $t\neq 0$ are essential.

Following the scheme of Section~2, we first   find a link between the discrete \e{eq:B1h} and continuous \e{eq:B6h}    realizations of Hankel operators. 
 Then we consider the case
    where    the matrix elements and the kernels of   Hankel operators admit the integral representations
 \begin{equation}
a_{n } = \int_{\cl{\Bbb D}}z^n dM(z), \q n=0,1,\ldots ,
 \label{eq:CM1}\end{equation}
 and
   \begin{equation}
{\bf a}(t) = \int_{\cl{{\Bbb  C}^+}} e^{-\zeta t} d\Sigma(\zeta), \q  t>0,
 \label{eq:CM2}\end{equation}
 with some complex measures $dM(z)$ and $ d\Sigma(\zeta)$.
Here $\Bbb D=\{z\in \Bbb C : |z| <1\}$ is the unit disc, ${\Bbb C}^+=\{  \zeta\in \Bbb C : \Re \zeta>0\}$ is the right half-plane, and $\cl{\Bbb D}$, $\cl{{\Bbb  C}^+}$ are the closures of these sets. Formulas \e{eq:CM1} and \e{eq:CM2} unify different types of integral representations of $a_{n } $ and ${\bf a}(t)$, for example, representations in terms of Carleson measures or in terms of symbols of the corresponding bounded Hankel operators.

 The central result of Section~3, Theorem~\ref{HQ}, formally means  that the ``operators" $H$ and $\bf H$ are unitarily equivalent provided the  measures $dM(z)$ and $ d\Sigma(\zeta)$ in \e{eq:CM1} and \e{eq:CM2} are linked by the equality
  \begin{equation}
d M(z)= 2\alpha  (\zeta +\alpha)^{-2} d\Sigma (\zeta), \q z=  \frac{\zeta-\alpha}{\zeta+\alpha},
 \label{eq:CM3}\end{equation}
 for some value of the parameter $\alpha>0$. Although quite simple, Theorem~\ref{HQ} is very  useful because it  relates the discrete and continuous representations directly avoiding the general construction of Section~2.
 More important, it allows one to translate spectral results obtained for the operator $\bf H$ into the results  for the operator $ H$, and vice  versa. Such examples are discussed in Section~4.
 
 \medskip

{\bf 1.3.}
 Section~4 is devoted to Hankel operators generalizing in different directions  two classical examples: the Hilbert matrix and the Carleman operator. To put our results into the right context, let us briefly recall basic spectral properties of these operators.  The Hilbert matrix is the Hankel operator $H$ defined by formula \e{eq:B1h} where
 $  a_n=  (n+1)^{-1}$ for all $n\geq 0$.  As shown in the papers \cite{Ma, Ro}, the spectrum of   $H$ is absolutely continuous,  it is simple  and coincides with the interval $[0,\pi]$.
   The Carleman operator is defined by formula \e{eq:B6h} where ${\bf a}(t)=t^{-1}$.  Using the Mellin transform, it is easy to show that  the spectrum of the operator $\bf H$ is absolutely continuous, has multiplicity $2$, and it also coincides with the interval $[0,\pi]$. So both these operators are bounded but not compact. It can be deduced from the results on the Hilbert matrix that Hankel operators $H$ are bounded if $   a_n =O (n^{-1})$ and they are compact if $   a_n =o (n^{-1})$   as  $   n\to \infty$. Similarly, the results on the
  Carleman operator imply  that Hankel operators $\bf H$ with integral kernels  ${\bf a}\in L^\infty_{\rm loc} ({\Bbb R})$ are bounded  if
  ${\bf a} (t)=O(t^{-1})$ and they are compact if
  ${\bf a} (t)=o(t^{-1})$ as $t\to 0$ and $t\to\infty$.

 In Section~4 we study   Hankel operators   $H$  with matrix elements $a_{n}$ such that  $a_{n} n\to\infty$.
 These operators are unbounded.
  We exhibit two quite different cases where the spectral analysis of Hankel operators   $H$  can be carried out sufficiently explicitly.  Our approach relies on Theorem~\ref{HQ} and the results on Hankel operators   $\bf H$ with singular integral kernels  ${\bf a}(t)$ obtained earlier in \cite{Yf1a, Ya}.
  
 We  emphasize that some properties of Hankel operators are more transparent in the discrete representation while 
 other properties -- in the continuous representation. Such examples are given in Section~4. So when studying 
 Hankel operators, it is very useful to keep in mind their various  representations.

 \section{Various  representations of convolutions and Toeplitz operators}

 {\bf 2.1.}
 First, we recall standard relations between various spaces we consider. Let us introduce the following diagrams: 
 \begin{equation}
  \begin{CD}
   f=\{f_{n}\}_{n\in{\Bbb Z}}   @ >>>  {\bf f} (x)=(\Phi  {\bf u})(x)     
   \\
    @VV    V     @   AA  A   
   \\
 u(\mu)    = ({\cal F}^* f) (\mu)  @>>>    {\bf u}(\lambda)= ({\cal U} u )(\lambda)
    \end{CD}
        \quad\q \q
        \begin{CD}
         \ell^2  ({\Bbb Z })   @> {\cal L}>>  L^2 ({\Bbb R};dx)
          \\
    @VV {\cal F}^*  V     @   AA \Phi A 
         \\
 L^2  ({\Bbb T})    @> {\cal U}>>   L^2  ({\Bbb R};d\lambda)    
    \end{CD}  
        \label{eq:dia}\end{equation}
 Here the unitary mapping ${\cal F}: L^2 ({\Bbb T})\to \ell^2({\Bbb Z})$ corresponds 
to expanding  a function in the Fourier series:
\[
 f_{n} =({\cal F} u)_{n} : =  \int_{\Bbb T} u(\mu)    \mu^{-n } d {\bf m}_{0}(\mu)
\]
where 
$$
d{\bf m}_{0}(\mu) = (2\pi i  \mu)^{-1}d\mu
$$
 is  the  normalized Lebesgue measure on the unit circle ${\Bbb T}$. The adjoint operator ${\cal F}^* \colon \ell^2 ({\Bbb Z}) \to L^2 ({\Bbb T})$ acts by the formula
   \begin{equation}
u(\mu) = ({\cal F}^* f)(\mu) =\sum_{n\in {\Bbb Z}}   f_{n}  \mu^n.
\label{eq:B3}\end{equation}
Similarly,  $\Phi $ is the Fourier transform, 
\[
{\bf f} (x)  = (\Phi {\bf u}) (x) : =(2\pi)^{-1/2}\int_{-\infty}^\infty e^{-ix\lambda}  {\bf u} (\lambda) d\lambda.
  \]
  Of course the operator $\Phi: L^2 ({\Bbb R; d\lambda})\to L^2 ({\Bbb R; dx})$ is unitary.
  The unitary operator  ${\cal U}={\cal U}_{\alpha}: L^2({\Bbb T})\to  L^2 ({\Bbb R}; d\lambda)$  is defined by the equality
\begin{equation}
({\cal U} u)(\lambda)= \sqrt{\tfrac{\alpha}{\pi}}  (\lambda+i \alpha)^{-1} u \bigl(\tfrac{\lambda-i\alpha}{\lambda+i\alpha}\bigr)
\label{eq:pwtq}\end{equation}
where a positive parameter $\alpha$ can be fixed in an arbitrary way.

Let us set ${\cal L}=\Phi {\cal U} {\cal F}^*   $. Then ${\cal L}  \colon \ell^2 ({\Bbb Z})\to \ell^2({\Bbb R}; dx)$ is
the unitary operator, and it    can be expressed    in terms of    the Laguerre polynomials.  Recall that  the Laguerre polynomials (see the book \cite{BE}, Chapter~10.12) are  defined by the formula
\[
{\sf L}_{n}^p  (t)= n!^{-1} e^t  t^{-p} d^n (e^{-t} t^{n +p  }) / dt^n , \q n=0,1, \ldots, \q t\geq 0, \q p>-1.
\]
Of course the polynomial ${\sf L}_{n}^p  (t)$ has  degree $n$; in particular, ${\sf L}_{0}^p  (t)=1$. These polynomials are orthogonal with respect to the measure $t^{p}e^{-t} dt$ and
\begin{equation}
\int_0^{\infty}{\sf L}_{n}^p  (t) {\sf L}_{m}^p  (t)t^{p}e^{-t} dt
= \frac{\Gamma (n+p+1)}{n!} \d_{n,m}
\label{eq:norm}\end{equation}
where $\Gamma (\cdot)$ is the gamma function and $\d_{n,m}$ is the Kronecker symbol.
The parameter $p>-1$ is  arbitrary, but we need the cases $p=0$ and $p =1$ only.

  Let us use the identity (see formula (10.12.32) in \cite{BE}) for ${\sf L}_{n}:={\sf L}_{n}^0$:
 \begin{equation}
  \int_{0}^\infty     {\sf L}_{n} (t ) e^{-(1/2+\zeta)t} dt=  \frac{1}{\zeta+1/2}\Big( \frac{2 \zeta-1}{2\zeta+1} \Big)^n, \q \Re \zeta> -1/2 .
\label{eq:KL4x}\end{equation}
Putting here $\z= -i\lambda$ and making the inverse Fourier transform, we find that
 \begin{equation}
 {\sf L}_{n} (2 \alpha x ) e^{-\alpha x} \1_{+}(x)= i(2\pi)^{-1}   \int_{-\infty}^\infty e^{- i x \lambda} (\lambda+i \alpha)^{-1}   \bigl(\frac{\lambda-i\alpha}{\lambda+i\alpha}\bigr)^n d\lambda, \q n=0,1, \ldots,
\label{eq:KL4x1}\end{equation}
where $\1_{+}(x)$ is the characteristic function of ${\Bbb R}_{+}$.
  Recall that  the Hardy space $ {\Bbb H}^2_{+}  ({\Bbb T})  $ (resp. $ {\Bbb H}^2_{-}  ({\Bbb T})  $) consists of functions $u\in  L^2   ({\Bbb T})  $ whose Fourier coefficients $({\cal F}u)_{n}=0$ for $n<0$ (resp., for $n\geq 0$).
Since the functions $\mu^n$, $ n =0,1,\ldots$,  form an orthonormal basis in the space ${\Bbb H}_{+}^2 ({\Bbb T})$, it follows from relations  \e{eq:pwtq} and \e{eq:KL4x1} that the functions $ -i (2\alpha)^{1/2} {\sf L}_{n} (2\alpha t ) e^{-\alpha t}  $  is an orthonormal basis in the space $L^2 ({\Bbb R}_{+})$.  We also see that
  $ i (2\alpha)^{1/2} {\sf L}_{n} (-2\alpha t ) e^{\alpha t}  $, $n=0,1,2, \ldots$,   is an orthonormal basis in the space $L^2 ({\Bbb R}_{-})$. Moreover,  relations  \e{eq:pwtq} and \e{eq:KL4x1} imply that the   operator ${\cal L} \colon \ell^2({\Bbb Z} )\to  L^2({\Bbb R} ;dx)$ acts by the formula  
\begin{multline}
({\cal L}  f ) (x)=  -i (2\alpha)^{1/2}  \sum_{n=0}^\infty f_{n}  {\sf L}_{n}  (2\alpha x) e^{-\alpha x}  \1_{+}(x)
\\
+  i (2\alpha)^{1/2}  \sum_{n= 0}^{\infty} f_{-n-1}  {\sf L}_{n}  (-2\alpha x) e^{\alpha x} \1_{+}(-x) , \q  f=\{ f_{n}\}_{n\in{\Bbb Z}}.
\label{eq:K3}\end{multline}
To be precise, the unitary operator ${\cal L} $ is first defined on the dense set ${\cal D}\subset \ell^2({\Bbb Z} )$ consisting of elements $f$ with only a finite number of non-zero components $f_{n}$, and then it is extended by the continuity onto the whole space $\ell^2({\Bbb Z} )$.

          \medskip

 {\bf 2.2.}
 Next, we discuss   representations of the  convolution/multiplication operators   in all these spaces. We start with the space  $\ell^2  ({\Bbb Z})$ where the operator of the discrete  convolution acts by the formula
  \begin{equation}
(Af )_{n}=\sum_{m\in {\Bbb Z}} a_{n-m} f_{m}, \q  f=\{ f_{n}\}_{n\in{\Bbb Z}}.
\label{eq:B1}\end{equation}
Without some assumptions on the sequence $a=\{ a_{n}\}_{n\in{\Bbb Z}}$, in general $Af\not\in \ell^2  ({\Bbb Z})$ even for $f\in  {\cal D}$.  Therefore instead of the operator $A $, we consider its quadratic form
  \begin{equation}
 a [f, f]= \sum_{n,m  \in{\Bbb Z}} a_{n-m}    f_{m } \ov{ f_{n}},\q f\in  {\cal D},
\label{eq:B2}\end{equation} 
that  consists of  a finite number of terms  for an arbitrary  sequence $a$.

Similarly,  the   convolution operator $ {\bf A}  $ acts in the space $L^2 ({\Bbb R}; dx)$ by the formula 
   \begin{equation}
({\bf A} {\bf f} )(x) = \int_{\Bbb R} {\bf a}(x-y) {\bf f}(y)  dy ,
\label{eq:B6}\end{equation}
and its quadratic form is given by the equality
 \begin{equation}
 {\bf a}[ {\bf f},  {\bf f}]  = \int_{\Bbb R} {\bf a}(x )  \Big(\int_{\Bbb R}       {\bf f}( y) \ov{ {\bf f}(x+y)} dy\Big) dx  
\label{eq:B7}\end{equation}
for  test functions 
 $ {\bf f}\in C_{0}^\infty ({\Bbb R})$.
Since, for such functions  $ {\bf f}$, the   function
 \[
 {\bf F}( x)= \int_{\Bbb R}       {\bf f}( y) \ov{ {\bf f}(x+y)} dy
 \]
 also belongs to $   C_{0}^\infty ({\Bbb R})$,
the form \e{eq:B7} is correctly defined for  a  distribution ${\bf a} $ in the space
$C_{0}^\infty ({\Bbb R})' $ dual to $C_{0}^\infty ({\Bbb R})$.
 
 Of course  the Fourier transform  allows one to realize convolutions as multiplication operators.  
Let $ {\cal P}= {\cal F}^* {\cal D}$ be the set  of all  quasi-polynomials \e{eq:B3}. For  a distribution $b\in  {\cal P}'$ (the space dual to $\cal P$), we formally define  the  operator $B $   in the space $L^2 ({\Bbb T})$
 by the equality
  \begin{equation}
 (B u) (\mu)=  b(\mu)u (\mu) ,
\label{eq:B4b}\end{equation} 
or, in precise terms, we introduce  its
 quadratic form   
   \begin{equation}
 b [u,u]= \int_{\Bbb T} b(\mu)|u(\mu)|^2 d{\bf m}_{0}(\mu), \q u\in {\cal P}. 
\label{eq:B4}\end{equation} 
Then   $B={\cal F}^* A {\cal F}$ if 
   \begin{equation}
b(\mu) = ({\cal F}^* a)(\mu) =\sum_{n\in {\Bbb Z}}   a_{n}  \mu^n , \q  a=\{ a_{n}\}_{n\in{\Bbb Z}}.
\label{eq:BX}\end{equation}
Strictly speaking, we have a relation between the  quadratic forms
\[
   b[u,u] =a[f,f] \q {\rm if} \q  f ={\cal F} u.
\]

Similarly, we put ${\cal Z}    := \Phi^* C_{0}^\infty ({\Bbb R})$. Recall that the set ${\cal Z}   $ consists of analytic functions satisfying a certain estimate at infinity (see, e.g., \cite{GUEVI+} for details).
For  a distribution $ {\bf b}\in  {\cal Z}'$ (the space dual to $\cal Z$), we formally define  the multiplication  operator $ \bf B $   in the space $L^2 ({\Bbb R}; d\lambda)$
 by the equality
 \begin{equation}
({\bf B} {\bf u})(\lambda) ={\bf b}(\lambda){\bf u} (\lambda) , 
\label{eq:B5b}\end{equation}
or, in precise terms, we introduce  its
 quadratic form   
 \begin{equation}
{\bf b} [{\bf u},{\bf u}]= \int_{\Bbb R} {\bf b}(\lambda)|{\bf u} (\lambda)|^2 d\lambda. 
\label{eq:B5}\end{equation} 
Then $ {\bf B} =\Phi^* {\bf A} \Phi$ if $   {\bf b} = (2\pi)^{1/2}\Phi^* {\bf a}$ where ${\bf a} \in C_{0}^\infty ({\Bbb R})'$. Strictly speaking, we have a relation between the  quadratic forms
\[
{\bf b} [{\bf u},{\bf u}] ={\bf a} [{\bf f},{\bf f}]    \q {\rm if} \q  {\bf f} =\Phi  {\bf u}.
\]

 Obviously, the operators $A$, $\bf A$, $B$ and $\bf B$ are formally symmetric if $a_{-n}=\ov{a_{n}}$, ${\bf a} (-x)=\ov{{\bf a} (x)}$,  $b(\mu)=\ov{b(\mu)}$ and ${\bf b}(\lambda)=\ov{{\bf b}(\lambda)}$. To be more precise, this means that the corresponding quadratic forms are real.

We emphasize that the bases $u_{n} (\mu)=\mu^n$, $n\in{\Bbb Z}$, and ${\cal F}u_{n}  $ are  the canonical bases in the spaces   $L^2 ({\Bbb T})$ and $\ell^2  ({\Bbb Z})$, respectively. On the contrary, the bases ${\bf u}_{n}={\cal U}u_{n}  $   and,  especially,   
$\Phi{\bf u}_{n}  $ in the spaces   $L^2 ({\Bbb R}; d\lambda)$ and  $L^2 ({\Bbb R}; d x)$ do not apparently play any distinguished role. So, it seems more natural to consider the form defined by  \e{eq:B7} on functions $ {\bf f}\in C_{0}^\infty ({\Bbb R})$ for  a distribution  ${\bf a}\in C_{0}^\infty ({\Bbb R})'$. Similarly, we consider the form  \e{eq:B5} for   $ {\bf u}\in {\cal Z} $ and ${\bf b}\in {\cal Z}'$.

 \medskip

 {\bf 2.3.}
Let us   find a link between the discrete and continuous  representations.  It is formally quite simple for the operators
$B$ and ${\bf B}$.
    By  definitions \e{eq:pwtq}  and  \e{eq:B4b}, the operator ${\bf B} ={\cal U}B {\cal U}^*$ acts   in the space $L^2 ({\Bbb R}; d\lambda)$
  as the multiplication by the function (distribution)
 \begin{equation}
{\bf b} (\lambda)= b \bigl(\tfrac{\lambda-i \alpha}{\lambda+i \alpha}\bigr)  .
\label{eq:B5bb}\end{equation}
Its quadratic form   is given by the formula \e{eq:B5}
   where
  \[
 {\bf u} (\lambda)= ({\cal U} u )(\lambda)= \sqrt{\tfrac{\alpha}{\pi}}  (\lambda+i\alpha)^{-1} \sum_{n\in {\Bbb Z}}f_{n} \bigl(\tfrac{\lambda-i\alpha}{\lambda+i\alpha}\bigr)^n 
\]
belongs to the set $\pmb{\cal P}= {\cal U}{\cal P}$ and   ${\bf b}$ belongs to the dual space
    $ \pmb{\cal P}'$. We emphasize however that this link is only formal because the sets $ {\cal Z}$ and $\pmb{\cal P}$ of test functions $ {\bf u} (\lambda)$ are different.
    
    It remains to directly link the representations of convolution operators in the spaces $\ell^2 ({\Bbb Z})$ and $L^2 ({\Bbb R};dx)$. Recall that the unitary operator   ${\cal L}=\Phi {\cal U} {\cal F}^* \colon \ell^2 ({\Bbb Z})\to L^2 ({\Bbb R};dx)$ is given by equality \e{eq:K3}. 
    Let  the operator $A $ be defined    in the space $\ell^2 ({\Bbb Z})$ by formula \e{eq:B1}  and ${\bf B}=  {\cal U} {\cal F}^* A {\cal F} {\cal U}^*$. Then $ {\bf A}= \Phi  {\bf B} \Phi^*  ={\cal L}A {\cal L}^*$ is  the   convolution   in the space $L^2 ({\Bbb R}; dx)$ acting by the formula  \e{eq:B6} where $ {\bf a}= (2\pi)^{-1/2}\Phi  {\bf b}$ and $ {\bf b}$ is defined by  \e{eq:BX},  \e{eq:B5bb}. The quadratic form of the operator $\bf A$ is defined by formula  \e{eq:B7}, but we have the same problem as for the multiplication operators: the domains $C_{0}^\infty ({\Bbb R})$ and ${\cal L}{\cal D}=\Phi \pmb{\cal P}$ of quadratic forms of the operators $ {\bf A}$ and ${\cal L}A {\cal L}^*$ are different.
         
    Our    goal is to find an expression for ${\bf a}(x)$, $x\in {\Bbb R}$, in terms of $a=\{ a_{n}\}_{n\in{\Bbb Z}}$. Recall 
    the relation   (see formula (10.12.15) in the book \cite{BE})
    \begin{equation}
 \frac{d}{dt}{\sf L}_{n} (  t ) =- {\sf L}_{n-1}^{1} ( t ) , \q t>0, \q \forall n\geq 1,
 \label{eq:Lag}\end{equation}
for the Laguerre polynomials.  Let ${\cal S}'$ be the   space   dual to the Schwartz space ${\cal S} = {\cal S} ({\Bbb R})$ of rapidly decaying $C^\infty$ functions.
 It follows from \e{eq:KL4x1} that
\begin{multline}
  \int_{-\infty}^\infty e^{- i x \lambda}    \bigl(\frac{\lambda-i\alpha}{\lambda+i\alpha}\bigr)^n d\lambda
  = i (\frac{d}{dx} +\alpha) \int_{-\infty}^\infty e^{- i x \lambda} (\lambda+i\alpha)^{-1}   \bigl(\frac{\lambda-i\alpha}{\lambda+i\alpha}\bigr)^n d\lambda
  \\
  = 
2 \pi   (\frac{d}{dx} +\alpha) \big(  {\sf L}_{n} (2\alpha x ) e^{-\alpha x} \1_{+}(x) \big)
  = 
-4  \pi  \alpha  {\sf L}_{n-1}^{1} (2\alpha x ) e^{-\alpha x} \1_{+}(x)  + 2\pi  \d(x)
\label{eq:B8}\end{multline}
where $\d(x)$ is the Dirac delta-function and the Fourier transform   is understood in the sense of ${\cal S}'$.
 Passing here to the complex conjugation and making the change of the variables $x \mapsto -x$, we also see that
\[
  \int_{-\infty}^\infty e^{- i x \lambda}    \bigl(\frac{\lambda-i\alpha}{\lambda+i\alpha}\bigr)^{-n} d\lambda
  = 
- 4  \pi  \alpha  {\sf L}_{n-1}^{1} (-2\alpha x ) e^{\alpha x} \1_{+}(-x)  + 2\pi  \d(x).
\]
Therefore it formally follows from  equalities  \e{eq:BX}  and \e{eq:B5bb}  that the distribution $   {\bf a} = (2\pi)^{-1/2}\Phi {\bf b}$ satisfies the relation
\begin{equation}
{\bf a}  (x)=\sum_{n\in{\Bbb Z}} a_{n} \d (x)
- 2\alpha   \sum_{n=1}^\infty    {\sf L}^{1}_{n-1}  (2\alpha | x|) e^{-\alpha |x|} \big( a_{n} \1_{+}(x) + a_{-n}   \1_{+}(-x) \big).
\label{eq:K3b}\end{equation}
An expression for this distribution can also be given in a  somewhat different form. Let $\varphi\in C_{0}^\infty({\Bbb R})$
(or ${\varphi\in\cal S}$). Then using \e{eq:Lag} and integrating by parts, we see that
\begin{multline}
\int_{-\infty}^\infty {\bf a}  (x) \varphi(x) dx=  a_{0} \varphi(0) 
\\
+   \sum_{n=1}^\infty  \int_{0}^\infty   {\sf L} _{n }  (2 \alpha  x) e^{-\alpha x} \big( a_{n} (\alpha \varphi(x) -  \varphi'(x))
+ a_{-n}  (\alpha \varphi(-x) +  \varphi'( - x)) \big) dx.
\label{eq:Kmod}\end{multline}
Equation   \e{eq:K3b} formally shows that
\[
{\bf a}  (x)=\kappa \d (x) + {\bf k}  (x)
\]
where
\begin{equation}
\sum_{n\in{\Bbb Z}} a_{n}=\kappa 
\label{eq:Kn2}\end{equation}
  and   ${\bf k} (x)$ is the second term in the right-hand side of \e{eq:K3b}.
   Using the identities \e{eq:norm}, we can solve equation \e{eq:K3b} for the coefficients $a_{n}$: 
\[
a_{\pm n}=- 2\alpha n^{-1}\int_{0}^{\infty}{\bf k}  (\pm x) {\sf L}^{1} _{n -1}  (2 \alpha  x) e^{-\alpha x} xdx, \q n=1,2, \ldots;
\]
then \e{eq:Kn2} yields $a_{0}$.

The relations between various representations of convolution/multiplication operators can be summarized by the following diagrams complementing \e{eq:dia}:
 \begin{equation}
    \begin{CD}
     a_{n}    @>>>   {\bf a} (x) 
     \\
        @VV    V     @   AA A    
     \\
   b (\mu)   @>>> {\bf b}(\lambda) 
 \end{CD}
  \q\q \q
      \begin{CD}
      A   @>>>  {\bf  A }  =\Phi  {\bf  B } \Phi^*
      \\
       @VV    V     @   AA   A     
   \\
B ={\cal F}^*  A  {\cal F}  @>>> {\bf  B } = {\cal U}  B {\cal U}^*
   \end{CD}
     \label{eq:DIA}\end{equation}

 \medskip
 
  {\bf 2.4.}
 Let us now consider Wiener-Hopf (Toeplitz) operators.  We start with the space  $\ell^2  ({\Bbb Z_{+}})$ where Wiener-Hopf operators $W$ act    again by  formula  \e{eq:B1} but now $n,m\in{\Bbb Z}_{+}$:
  \begin{equation}
(W f )_{n}=\sum_{m\in {\Bbb Z}_{+}} a_{n-m} f_{m}, \q  f=\{ f_{n}\}_{n\in{\Bbb Z}_{+}}.
\label{eq:B1w}\end{equation}
Quite similarly to \e{eq:B2}, their quadratic forms are defined by the formula
  \begin{equation}
 w [f, f]=     \sum_{n, m\in{\Bbb Z}_{+}} a_{n-m} f_{m} \ov{ f_{n}}
\label{eq:B2w}\end{equation} 
where $f\in  {\cal D}_{+} : = {\cal D}  \cap  \ell^2({\Bbb Z}_{+} )$ and the sequence $a=\{ a_{n}\}_{n\in{\Bbb Z}}$ is again arbitrary.

A  Wiener-Hopf  operator $ {\bf W}  $  acts in the space $L^2 ({\Bbb R}_{+} )$ by the formula 
   \[
({\bf W} {\bf f} )(x) = \int_{{\Bbb R}_{+}} {\bf a}(x-y) {\bf f}(y)  dy ,
\]
and its quadratic form is again given by the equality
 \begin{equation}
 {\bf w}[ {\bf f},  {\bf f}]  =  
   \int_{\Bbb R} {\bf a}(x )  \Big(\int_{\Bbb R}       {\bf f}( y) \ov{ {\bf f}(x+y)} dy\Big) dx
\label{eq:B7W}\end{equation}
 where   ${\bf a}\in C_{0}^\infty ({\Bbb R})'$ but now
 $ {\bf f}\in C_{0}^\infty ({\Bbb R}_{+})$.

Next, we pass to the representation of Toeplitz operators in the Hardy spaces. 
Denote by $P_{+}$ the orthogonal projection in $ L^2   ({\Bbb T})  $ onto the Hardy space $ {\Bbb H}^2_+  ({\Bbb T})  $, and  let, as before, the operator  $B$ be formally defined by equality \e{eq:B4b} where the  distribution $b\in {\cal P}'$.    Then  the Toeplitz operator $T: {\Bbb H}^2_{+}  ({\Bbb T})  \to {\Bbb H}^2_{+}  ({\Bbb T})  $ is defined by the relation 
 \[
 T u =P_{+} B u 
 \]
 on elements  $ u\in {\cal P}_{+}:= {\cal F}^* {\cal D}_{+}$; obviously, the set ${\cal P}_{+}$ consists   of all  polynomials $u(\mu)=\sum_{n\in {\Bbb Z}_{+}} f_{n} \mu^n $. 
 The  quadratic form $ t [u,u]$ of the  operator $T$ is given by the right-hand side of \e{eq:B4} where $ u\in {\cal P}_{+}$ and $b\in {\cal P}'$ are arbitrary.

 Finally, we discuss the  representation in the Hardy space $ {\Bbb H}^2_+  ({\Bbb R})  $ consisting of functions 
  $  {\bf u}\in L^2   ({\Bbb R} )  $ whose Fourier transforms $ (\Phi {\bf u})(x)=0$ for $x<0$.  Let  ${\bf P}_+$ be the orthogonal projection in $ L^2   ({\Bbb R} )  $ onto   $ {\Bbb H}^2_+  ({\Bbb R})  $, and let ${\bf B}$ be the   operator \e{eq:B5b}.
  Then  the Toeplitz operator ${\bf T}: {\Bbb H}^2_{+}  ({\Bbb R})  \to {\Bbb H}^2_{+}  ({\Bbb T})  $ is defined by the relation 
 \[
{\bf T} {\bf u}  = {\bf P}_{+} {\bf B} {\bf u}  .
 \]
Its quadratic form ${\bf t} [{\bf u} ,{\bf u} ]$ is given by the the right-hand side of   
\e{eq:B5} where $ {\bf u} \in   \Phi^* C_{0}^\infty ({\Bbb R}_{+})=: {\cal Z}_{+}$ and ${\bf b}\in  {\cal Z}  '$.   The relation \e{eq:B5bb} between the functions (distributions) $b(\mu)$ and ${\bf b} (\lambda)$ remains of  course true.

Note that the Laguerre operator ${\cal L} \colon \ell^2({\Bbb Z}_\pm )\to  L^2({\Bbb R}_\pm  )$ (here ${\Bbb Z}_-
={\Bbb Z}\setminus {\Bbb Z}_{+}$) and
\begin{equation}
({\cal L}  f ) (t)= - i (2\alpha)^{1/2}  \sum_{n=0}^\infty f_{n}  {\sf L}_{n}  (2\alpha t) e^{-\alpha t}   , \q  f=\{ f_{n}\}_{n\in{\Bbb Z}_{+}}, \q t>0,
\label{eq:K3TH}\end{equation}
but formula  \e{eq:K3b} for kernel of the operator $ {\bf  W }$ remains true.

For  Toeplitz operators,  instead of  \e{eq:dia}, \e{eq:DIA} we have  the  diagrams 
 \begin{equation}
       \begin{CD}
         \ell^2  ({\Bbb Z }_+)   @> {\cal L}>>  L^2 ({\Bbb R}_{+})
          \\
    @VV {\cal F}^*  V     @   AA \Phi A 
         \\
 {\Bbb H}_{+}^2  ({\Bbb T})    @> {\cal U}>>    {\Bbb H}_{+}^2  ({\Bbb R})  
    \end{CD}
    \q\q\q
    \begin{CD}
      W    @>>>  {\bf  W }  =\Phi  {\bf  T }  \Phi^*
      \\
       @VV    V     @   AA A    
   \\
T={\cal F}^*  W  {\cal F}   @>>> {\bf  T } = {\cal U}  T {\cal U}^*  
    \end{CD}
      \label{eq:DIAT}\end{equation}
 
   Of course all the remarks above concerning a certain difference between
     the discrete and continuous representations of  convolution operators apply also to Wiener-Hopf  operators.   The case of semibounded Wiener-Hopf  operators is specially discussed
     in \cite{YaT} (the discrete  representation) and in \cite{YaWH} (the continuous  representation).

    \medskip
 
  {\bf 2.5.}
  Let us say a few words about bounded operators. For Wiener-Hopf operators $W$ defined via the quadratic form  \e{eq:B2w}, it is  the classical Toeplitz result  that the operator $W$ is bounded, that is,
  \[
 | w[f,f] | \leq C \| f\|^2
 \]
for some $C>0$, if and only if $a={\cal F} b$ where $b\in L^\infty ({\Bbb T})$.  The corresponding result for integral Wiener-Hopf operators $\bf W$ is stated explicitly in \cite{YaWH}. In the assertion below the Fourier transform   is understood in the sense of the Schwartz space ${\cal S}'$.

   \begin{proposition}\label{B-H}\cite[Theorem~1.1]{YaWH} 
  Let the form ${\bf w}[{\bf f}, {\bf f}] $ be defined by the relation \e{eq:B7W} where $ {\bf f}\in C_{0}^\infty({\Bbb R}_{+})$ and the distribution $ {\bf a}\in C_{0}^\infty({\Bbb R})'$. Then
    \[
 | {\bf w}[{\bf f}, {\bf f}] | \leq C \| {\bf f} \|^2,
 \] 
 $($so  that the corresponding operator $\bf W$ is bounded$\,)$ if and only if $   {\bf a} = (2\pi)^{-1/2}\Phi {\bf b}$ where ${\bf b}\in L^\infty ({\Bbb R} )$. Moreover, $\|\bf W\|= \| {\bf b}\|_{L^\infty ({\Bbb R} )}$.
      \end{proposition}

 Observe that Proposition~\ref{B-H} is not a direct consequence  of  the   Toeplitz criterion for the boundedness of the operators $W$ in the space $\ell^{2} ({\Bbb Z}_{+})$. The difference is that the domains ${\cal D}$ and $C_{0}^\infty({\Bbb R}_{+})$    of the corresponding quadratic forms are not linked by the Laguerre transform ${\cal L}$.

\section{ Hankel operators}

 {\bf 3.1.}
 Let us pass to Hankel operators.   We start with the space  $\ell^2  ({\Bbb Z_{+}})$ where Hankel  operators $H$ act    according to the formula (cf. \e{eq:B1w})
  \[
(H f )_{n}=\sum_{m\in {\Bbb Z}_{+}} a_{n+m} f_{m}, \q  f=\{ f_{n}\}_{n\in{\Bbb Z}_{+}},
\]
  and their quadratic forms (cf. \e{eq:B2w}) are defined by the formula
  \begin{equation}
 h [f, f]=     \sum_{n, m\in{\Bbb Z}_{+}} a_{n+m} f_{m} \ov{ f_{n}}
\label{eq:B2h}\end{equation} 
where $f\in  {\cal D}_{+}= {\cal D}  \cap  \ell^2({\Bbb Z}_{+} )$ and the sequence $a=\{ a_{n}\}_{n\in{\Bbb Z}_{+}}$ is  arbitrary.

A  Hankel operator $ {\bf H}  $ acts in the space $L^2 ({\Bbb R}_{+})$ by the formula 
    \[
({\bf H} {\bf f} )(t) = \int_{{\Bbb R}_{+}} {\bf a}(t+s) {\bf f}(s)  ds ,
\]
and its quadratic form is given by the equality
 \begin{equation}
 {\bf h}[ {\bf f},  {\bf f}]  =  \int_{{\Bbb R}_{+}} {\bf a}(t )  \Big(\int_{0}^{t}       {\bf f}( s) \ov{ {\bf f}(t-s)} ds\Big) dt .
\label{eq:B7h}\end{equation} 
  Since, for all test functions 
 $ {\bf f}\in C_{0}^\infty ({\Bbb R}_{+})$, the   function
  \begin{equation}
{\bf F}( t) = \int_{0}^{t}       {\bf f}( s) \ov{ {\bf f}(t-s)} ds
\label{eq:B7hh}\end{equation} 
also belongs to $   C_{0}^\infty ({\Bbb R}_{+})$,
the form \e{eq:B7h} is correctly defined for all   distributions ${\bf a}\in C_{0}^\infty ({\Bbb R}_{+})'$ in the space dual to $C_{0}^\infty ({\Bbb R}_{+})$.

   Let us now pass to the representation of Hankel operators in the  Hardy spaces 
${\Bbb H}^2_{+} ({\Bbb T}) $ and ${\Bbb H}^2_{+} ({\Bbb R}) $.   A Hankel  operator $ G$ in the space ${\Bbb H}^2_{+} ({\Bbb T})$  is formally defined by the relation 
 \[
 G u =P_{+} B Ju 
 \]
 where $(J u) (\mu) =\bar{\mu} u(\bar{\mu})$ so that the involution $J\colon {\Bbb H}^2_{\pm} ({\Bbb T})\to {\Bbb H}^2_{\mp} ({\Bbb T})$. As before   $B $ is  the  multiplication  operator  \e{eq:B4b}  by the function $b(\mu) $ in the space $L^2 ({\Bbb T})$. To be precise, we define $G$ via its  quadratic form
 \begin{equation}
 g[ u,u]= \int_{\Bbb T} \omega (\mu) u(\bar{\mu}) \ov{u(\mu)}    d{\bf m}_{0}(\mu) 
\label{eq:QF}\end{equation}
 where  $   \omega(\mu)=\bar{\mu}b(\mu)$,   $ u\in   {\cal P}_{+} =  {\cal F}^{*} {\cal D}_{+}$ and $\omega \in  {\cal P}_{+}'$. Obviously, we have
  \begin{equation}
g[   u,  u]=  h[ {\cal F}u, {\cal F}u] \q {\rm if} \q  \omega= {\cal F}^{*} a.
\label{eq:QF1}\end{equation}

   Hankel operators in the  Hardy space 
${\Bbb H}^2_{+} ({\Bbb R})\subset L^2 ({\Bbb R})$ 
of functions analytic in the upper half-plane are defined quite similarly. Let, as before,  ${\bf P}_+$ be the orthogonal projection in $ L^2   ({\Bbb R} )  $ onto   $ {\Bbb H}^2_+ ({\Bbb R})  $, and let ${\bf B}$ be the   operator \e{eq:B5b}. A Hankel  operator $\bf G$   is formally defined  in the space ${\Bbb H}^2_{+} ({\Bbb R})$ by the relation 
 \[
 {\bf G} {\bf u} ={\bf P}_{+} {\bf B} {\bf J} {\bf u} 
 \]
 where $(  {\bf J} {\bf u} ) (\lambda) ={\bf u}(-\lambda)$. To be precise, we have to pass to the Fourier transform  in \e{eq:B7h}  which yields the representation
\begin{equation}
 {\bf g }[{\bf u}, {\bf u}]= \int_{\Bbb R}\Omega  (\lambda) {\bf u} (-\lambda)   \ov{{\bf u} (\lambda)}    d\lambda, \q
 {\bf u}\in {\cal Z}_{+}  ,  
\label{eq:QFR}\end{equation}
for the quadratic form of the operator $\bf G$. Obviously, we have
  \begin{equation}
   {\bf g }[{\bf u}, {\bf u}]=  {\bf h }[ \Phi{\bf u}, \Phi{\bf u}]
 \q {\rm if} \q      \Omega =(2\pi)^{1/2}\Phi^{*} {\bf a} .
\label{eq:QF1c}\end{equation}

The functions $\omega  (\mu)$ and $\Omega  (\lambda)$ are known as symbols of the discrete and continuous Hankel operators. Making the change of the variables $\mu=(\lambda-i\alpha) (\lambda+i\alpha)^{-1}$ in \e{eq:QF} we see that
\[
 {\bf g }[{\cal U}u, {\cal U}u]=  g[ u,u]
\]
if
 \begin{equation}
\Omega (\lambda)= - \bigl(\tfrac{\lambda-i \alpha}{\lambda+i \alpha}\bigr) \omega \bigl(\tfrac{\lambda-i \alpha}{\lambda+i \alpha}\bigr).  
\label{eq:B5h}\end{equation}
For Hankel operators $G$ and ${\bf G}$, this equality plays the role 
 of \e{eq:B5bb}. 
  Since   ${\bf a} = (2\pi)^{-1/2}\Phi \Omega$ and $\omega= {\cal F}^{*} a$, relation \e{eq:B8} where $x> 0$   yields the representation (cf. \e{eq:K3b})
\begin{equation}
{\bf a}  (t)= 
 2\alpha   \sum_{n=0}^\infty  a_{n}  {\sf L}^{1}_{n}  (2\alpha t) e^{-\alpha t}, \q t >0.
\label{eq:Khb}\end{equation}
Let us now use that $(n+1)^{-1/2} 2\alpha t^{1/2}  {\sf L}^{1}_{n}  (2\alpha t) e^{-\alpha t}$, $n=0,1,\ldots$, is the orthonormal basis in the space $L^{2}({\Bbb R}_{+})$. Therefore it follows from \e{eq:Khb} that
\begin{equation}
a_{n}   = 
 2\alpha   (n+1)^{-1}\int_{{\Bbb R}_{+}}  {\bf a}  (t) {\sf L}^{1}_{n}  (2\alpha t) e^{-\alpha t} tdt .
\label{eq:Khd}\end{equation}

  Similarly to  Toeplitz operators (cf. \e{eq:DIAT}),   the relations between various representations of Hankel operators can be summarized by the following diagram:
 \[
    \begin{CD}
     H    @>>>  {\bf  H }  =\Phi  {\bf  H }  \Phi^*
     \\
       @VV    V     @   AA A     
   \\
  G  ={\cal F}^*  H  {\cal F}  @>>> {\bf  G } = {\cal U}  G {\cal U}^*
    \end{CD}
    \]

It is easy to see that Hankel operators $H$, $\bf H$, $G$ and  $\bf G$  are formally symmetric  if  
          \[
          a_{n}=\bar{a}_{n} , \q  {\bf a}(t ) = \overline{{\bf a} (t )}, \q 
   \omega(\bar{\mu} )=\overline{\omega(\mu)} \q    \q{\rm and }\q \Omega (-\lambda)=\overline{\Omega(\lambda)}.
    \]
    To be more precise, this means that the corresponding quadratic forms are real.

 Similarly to the case of convolutions and Toeplitz operators, the relations \e{eq:Khb} and \e{eq:Khd} are only formal.
 We also note that the domains ${\cal D}_{+}$ and $C_{0}^{\infty} ({\Bbb R}_{+})$ of the forms \e{eq:B2h} and \e{eq:B7h} are not related by the Laguerre transform \e{eq:K3TH}; likewise, the domains $ {\cal P}_{+}$ and $
{\cal Z}_{+}$ of the forms \e{eq:QF} and\e{eq:QFR} are not related by the change of variables \e{eq:pwtq}. We emphasize however that this discrepancy is essential in singular cases only.

\medskip

{\bf 3.2.}
As far as the conditions of boundedness are concerned, the results on Hankel operators are formally similar to those on Toeplitz operators.  For Hankel operators $H$ defined via the quadratic form  \e{eq:B2h}, it is  the classical Nehari result  that the operator $H$ is bounded 
 if and only if there exists $\omega\in L^\infty ({\Bbb T})$ such that $a={\cal F} \omega$.  An important difference compared to Toeplitz operators is that now the equation  $a={\cal F} \omega$ does not determine the function $\omega$ uniquely and  may  also be satisfied with   unbounded   $\omega$.
 
The corresponding result for integral Hankel operators $\bf H$ is stated explicitly in \cite{Y}. 
  
   \begin{proposition}\label{NehC}\cite[Theorem~3.4]{Y}  
  Let the form ${\bf h}[ {\bf f} , {\bf f}] $ be defined by the relation \e{eq:B7h} where the distribution $ {\bf a}\in C_{0}^\infty({\Bbb R}_{+})'$. Then the corresponding operator $\bf H$ is bounded if and only if there exists a function $\Omega \in L^\infty ({\Bbb R})$ such that  
      \[
 {\bf a}= (2\pi)^{-1/2} \Phi \Omega.
\]
In this case $\|{\bf H}\|= \| \Omega\|_{L^\infty ({\Bbb R} )}$.
      \end{proposition}
      
      By the same reasons as for Toeplitz  operators (see the remark after Proposition~\ref{B-H}), this  result
        is not a direct consequence  of  the   Nehari theorem.  
              
      Of course, the operators $\cal F$, $\Phi$, $\cal U$   and $\cal L$ establishing the equivalence of various representations of Hankel operators are not unique. For example,  the operator ${\cal U}={\cal U}_{\alpha}$  defined by equality \e{eq:pwtq}  depends on the parameter $\alpha>0$, and there is no distinguished choice of this parameter. To state the problem precisely, for each of the spaces ${\cal H}=\ell^2 ({\Bbb Z}_{+}), L^2 ({\Bbb R}_{+}),  {\Bbb H}_{+}^2 ({\Bbb T}),   {\Bbb H}_{+}^2 ({\Bbb R})$, let us introduce the 
      group $ {\Bbb G}  ({\cal H} )$ of all unitary automorphisms of the set $ {\Bbb A}  ({\cal H} )$ of bounded Hankel operators in the space $\cal H$.
        By definition, a unitary operator $\sf U\in {\Bbb G} ({\cal H}) $ if and only if ${\sf U} {\sf H}{\sf U}^*  \in {\Bbb A}  ({\cal H} ) $    for all   operators ${\sf H} \in {\Bbb A} ({\cal H})$. For different $\cal H$, the groups ${\Bbb G} ({\cal H}) $ are related by  the transformations $\cal F$, $\Phi$, $\cal U$   and $\cal L$. So it suffices to describe this group for one of the choices of $\cal H$. 
        
        It turns out that the group ${\Bbb G} ({\Bbb H}_{+}^2 ({\Bbb R})) $ admits a very explicit description.  Let ${\sf D}_\rho$,
 \[
( {\sf D}_\rho {\bf u})(\lambda) = \rho^{1/2} {\bf u} (\rho \lambda), \q \rho>0,
\]
be the dilation operator in ${\Bbb H}_{+}^2 ({\Bbb R})$, and let the involution ${\sf J}  $   be defined in this space by the equation
 \[
( {\sf J} u)(\lambda) =i \lambda^{-1} u(-\lambda^{-1}).
\]
Obviously, ${\sf D}_\rho  {\bf G} {\sf D}_\rho ^* $ and $ {\sf J} {\bf G}  {\sf J}^* $  are Hankel operators for all Hankel operators ${\bf G}\in {\Bbb A} ({\Bbb H}_{+}^2 ({\Bbb R}))$. 
Surprisingly, the group $ {\Bbb G}  ({\Bbb H}_{+}^2 ({\Bbb R}))$ is exhausted by these transformations.
Let us state the precise result.

 \begin{theorem}\label{Aut}\cite[Theorem~A.1]{Yfr}
A unitary operator ${\sf U}\in  {\Bbb G}  ({\Bbb H}_{+}^2 ({\Bbb R}))$ if and only if it has one of the two forms:  ${\sf U}= \mu {\sf D}_\rho$ or ${\sf U}=\mu {\sf D}_\rho {\sf J} $ for some $ \mu \in {\Bbb T}$ and $\rho>0$.
   \end{theorem}

\medskip

{\bf 3.3.}
Now we consider the case when the matrix elements $a_{n}$ of a Hankel operator $H$ and the integral kernel ${\bf a} (t) $ of a Hankel operator $\bf H$ are given by formulas \e{eq:CM1} and \e{eq:CM2}, respectively.  Here $dM(z)$ is  a finite complex measure on   $\cl{\Bbb  D}$  and   $d\Sigma(\zeta)$ is a locally finite complex measure on   $\cl{{\Bbb C}^+}$. We suppose that the measures $dM(z)$    and   $d\Sigma(\zeta)$ are
linked by equation \e{eq:CM3}; in this case $M(\{1\})=0$.  We  denote by    $d |M|(z)$ and $d |\Sigma |(\zeta)$  the   variations of these measures and assume that 
 \begin{equation}
| M | (\cl{\Bbb  D})=2\alpha\int_{\cl{{\Bbb  C}^+}}|\zeta +\alpha|^{-2} d|\Sigma |(\zeta)<\infty.
 \label{eq:CM}\end{equation}
  Relation \e{eq:CM1}   only  implies  that the sequence $a_{n }$ is bounded as $n\to\infty$, and hence 
 the operator $H$ is not defined in $\ell^2 ({\Bbb Z}_{+})$ even on the set $\cal D   $. Similarly,  the operator $\bf H$ is not defined in $L^2({\Bbb R}_{+})$ even on the set $C_{0}^\infty({\Bbb R}_{+})$. So as usual, instead of operators we have to work with the corresponding quadratic forms \e{eq:B2h} and \e{eq:B7h}. 
 
 Formula \e{eq:CM2} defines ${\bf a}(t)$ as a distribution on a set of test functions, denoted $\cal X$, that can be chosen as follows. A function ${\bf F}\in{\cal X}$ if and only if ${\bf F}\in C^{\infty} ({\Bbb R}_{+})$, there exist limits ${\bf F}^{(k)}(+0)$ for $k=0,1,2$, ${\bf F}(+0)=0$ and ${\bf F}^{(k)}\in L^{1}({\Bbb R}_{+})$ for $k=0,1,2$. Integrating twice by parts, we find that, for such ${\bf F}$, the estimate
 \[
 \big|  \int_{0}^{\infty} e^{-\z t} {\bf F} (t) dt \big| \leq C |\zeta + 1 |^{-2},\q \z \in \cl{{\Bbb  C}^+},
 \]
 holds. Therefore the form
  \begin{equation}
\la {\bf a} , {\bf F} \ra :=\int_{\cl{{\Bbb  C}^+}} \big(  \int_{0}^{\infty} e^{-\z t} {\bf F} (t) dt \big) d\Sigma (\zeta) 
 \label{eq:CMF}\end{equation}
 is well defined for all ${\bf F}\in \cal X$.   
 
 Let us also introduce a set of test functions ${\bf f}\in C^{\infty} ({\Bbb R}_{+})$, denoted $\cal Y$, such that   there exist limits ${\bf f}^{(k)}(+0)$   and ${\bf f}^{(k)}\in L^{1}({\Bbb R}_{+})$ for $k=0,1,2$. By definition \e{eq:B7hh}, we have ${\bf F}\in{\cal X}$ if ${\bf f}\in{\cal Y}$. Therefore if  ${\bf h}\in {\cal X}'$,  then the quadratic form ${\bf h} [{\bf f} ,{\bf f} ]$  is well defined for all ${\bf f}\in \cal Y$.   Note that ${\bf D}_{+} : = {\cal L} {\cal D}_{+}\subset \cal Y$.

  Our main result in this subsection formally means  that the ``operators" $H$ and $\bf H$ are unitarily equivalent provided the corresponding measures are linked by the equality \e{eq:CM3}. 
 Let us state this result precisely.

 \begin{theorem}\label{HQ} 
Let  the matrix elements $a_n$ and the integral kernel ${\bf a}(t)$   be  given by equalities \e{eq:CM1} and \e{eq:CM2}, and let $h[f,f]$ and ${\bf h} [{\bf f} ,{\bf f} ]$ be the corresponding quadratic forms  \e{eq:B2h} and \e{eq:B7h}
defined for  $f\in {\cal D}_{+}$  and  ${\bf f}\in   {\cal Y}  $.
Suppose that the measures  $dM(z)$ and $d\Sigma (\zeta)$  are linked by the relation \e{eq:CM3} for some $\alpha >0$ and satisfy the condition  \e{eq:CM}. Then for all $f\in {\cal D}_{+}$ the identity
  \begin{equation}
h [f,f]= {\bf h} [{\cal  L }f , {\cal  L }f ]
\label{eq:ident}\end{equation}
holds.
    \end{theorem}
    
 \begin{pf}
 Let $f\in{\cal D}_{+}$ and ${\bf f} ={\cal L}f \in{\cal Y}$. It follows from formulas  \e{eq:B7hh}  and \e{eq:CMF}     that in this case   \begin{equation}
 {\bf h}[ {\bf f}, {\bf f}] =  \int_{\cl{{\Bbb  C}^+}}   d\Sigma(\zeta) \big(\int_{0}^\infty e^{-\zeta t} \ov{ {\bf f}(t)} dt \big)\big( \int_{0}^\infty e^{-\zeta t}  {\bf f}(s)  ds  \big).
 \label{eq:CM6}\end{equation} 
 Moreover, in view of the identity   \e{eq:KL4x} and definition   \e{eq:K3TH}, we have
   \[
 \int_{0}^\infty e^{-\zeta t}  {\bf f} (t) dt =   \frac{-i \sqrt{2\alpha}}{\zeta+\alpha} \sum_{n=0}^\infty f_{n}  \Big( \frac{ \zeta-\alpha}{\zeta+\alpha} \Big)^n
\]
 where the sum consists of a finite number of terms.
 Substituting this expression into relation \e{eq:CM6}, we find that
   \[
 {\bf h}[ {\bf f}, {\bf f}] =2\alpha \sum_{n,m=0}^\infty f_{m} \ov{f_{n}}\int_{\cl{{\Bbb  C}^+}}     
 \Big( \frac{ \zeta-\alpha}{\zeta+\alpha} \Big)^{n+m} \frac{1}{(\zeta+\alpha)^2}d\Sigma(\zeta).
\]
 After the change of the variables
 $z=\frac{\zeta-\alpha}{\zeta+\alpha}$, we see that this expression equals
  \[
 {\bf h}[ {\bf f}, {\bf f}] = \sum_{n,m=0}^\infty f_{m} \ov{f_{n}}\int_{\cl{\Bbb  D}}    
z^{n+m} dM(z)
\] 
 where the measure $dM(z)$ is defined by relation \e{eq:CM3}. According to \e{eq:CM1} this expression coincides with \e{eq:B2h}. This concludes the proof of the identity \e{eq:ident}.
   \end{pf}
   
   \begin{corollary}\label{HQb} 
   Suppose that the expression \e{eq:ident} is estimated by $C \| f\|^{2} $ with some $C>0$. Then there exist bounded operators
       $H$ and 
   $\bf H$ corresponding to these quadratic forms, and they are unitarily equivalent:
   ${\bf H} {\cal  L } = {\cal  L } H$.
    \end{corollary}

Obviously,  coefficients \e{eq:CM1} and   function   \e{eq:CM2} are real if the measures $dM(z)$ and $d\Sigma (\z)$ are invariant with respect to the complex conjugation:
\begin{equation}
 M(\ov{X})=\ov{M(X)}\q {\rm and} \q  \Sigma(\ov{Y})=\ov{\Sigma(Y)}
 \label{eq:sym}\end{equation}
 for all $X\subset\cl{\Bbb D}$ and $Y\subset\cl{\Bbb C}^{+}$.
 
  \begin{corollary}\label{HQc} 
  Under the assumptions of Theorem~\ref{HQ},
 let condition \e{eq:sym} be satisfied. Suppose that an operator $H$ is defined on ${\cal D}_{+}$ and $(Hf,f)=h[f,f]$ for 
 $f\in {\cal D}_{+}$, or equivalently an operator $\bf H$ is defined on ${\bf D}_{+}$ and $({\bf H}{\bf f}, {\bf f})= {\bf h}[ {\bf f},{\bf f}]$ for  ${\bf f}\in {\bf D}_{+}$. If $H$ is essentially self-adjoint  on ${\cal D}_{+}$ $($or equivalently  $\bf H$ is essentially self-adjoint  on ${\bf D}_{+})$, then their closures are unitarily equivalent:
   $(\cl{\bf H} )\, {\cal  L } = {\cal  L } (\cl{H})$.
    \end{corollary}

    In particular, the measures $dM(z)$ and $d\Sigma(\z)$ may be supported by the intervals $[ -1,1 ]$ and $[ 0,\infty)$, respectively; as before we suppose also that $M(\{1\}) =0$.     Then Theorem~\ref{HQ} reads as follows.

 \begin{theorem}\label{HQX} 
Let  the matrix elements $a_n$ and the kernel ${\bf a}(t)$   be  given by the equalities    
 \begin{equation}
a_{n } = \int_{-1}^1 \nu^n dM(\nu), \q n=0,1,\ldots, 
 \label{eq:W}\end{equation}
 and
      \begin{equation}
{\bf a}(t) = \int_{0}^\infty e^{-\lambda t} d\Sigma(\lambda), \q  t>0.
 \label{eq:S}\end{equation}
Suppose that the measures  $dM(\nu)$ and $d\Sigma (\lambda)$  are linked by the relation 
 \begin{equation}
d M(\nu)=2\alpha (\lambda + \alpha)^{-2} d\Sigma (\lambda), \q \nu=\frac{\lambda- \alpha}{\lambda+\alpha},
 \label{eq:MS}\end{equation}
  for some $\alpha >0$ and satisfy the condition  
  \[
| M | ([-1,1))=2\alpha\int_0^{\infty} | \lambda +\alpha|^{-2} d|\Sigma |(\lambda)<\infty.
\]
 Then for all $f\in {\cal D}_{+}$ the identity \e{eq:ident} holds.
    \end{theorem}
    
    The following result is a combination of Theorem~\ref{HQX}, of Theorems~1.2 and 3.4 in \cite{Yunb} on the discrete case and the preceding results (Theorem~3.10) of  \cite{Yf1a} on the continuous case.
    
   \begin{theorem}\label{HQXb} 
   Under the assumptions of Theorem~\ref{HQX}
   suppose that the measure \e{eq:MS} is non-negative and that $M(\{-1\})=\Sigma(\{0\})=0$. Then:
   \begin{enumerate}[\rm(i)]
   
   \item
   The form $h[f, f]$ defined on ${\cal D}_{+}$ is closable, and it
    is closed on the set of elements $f=(f_{0}, f_{1 }, \ldots)\in \ell^{2} ({\Bbb Z}_{+})$ such that
   \[
h[f,f]= \int_{-1}^{1} \Big|\sum_{n=0}^\infty  f_{n} \nu^n \Big|^2    dM(\nu) <\infty.
\]

   \item
   The form ${\bf h} [{\bf f}, {\bf f}]$ defined on $  {\bf D}_{+}$ is closable, and it
    is closed on the set of functions ${\bf f}\in L^{2} ({\Bbb R}_{+})$ such that
   \[
{\bf h} [{\bf f}, {\bf f}]= \int_0^\infty \Big| \int_0^\infty e^{-\lambda t} {\bf f}(t) dt\Big|^2    d\Sigma (\lambda) <\infty.
 \]

     \item
The non-negative operators $H$ and 
   $\bf H$ corresponding to these quadratic forms   are unitarily equivalent:
   ${\bf H} {\cal  L } = {\cal  L } H $.
    \end{enumerate}
    \end{theorem}

 Sometimes the representations \e{eq:CM1} and \e{eq:CM2} are  too restrictive. For example, if a Hankel operator $\bf H$ has kernel $a(t)=t^k e^{-\alpha t}$, $\Re\alpha>0$ (such $\bf H$ has a finite rank), then representation  \e{eq:CM2} is  formally satisfied with
 \[
 d\Sigma (\z)=\d^{(k)}(\z-\alpha)d{\bf M}_{0}(\z)
 \]
 ($d{\bf M}_{0}(\z)$ is the planar Lebesgue measure and $\d(\z)$ is the delta-function) which is a measure for $k=0$ only.
 A very general situation where ${\bf a}\in C_{0}^\infty
 ({\Bbb R}_{+})' $ is an arbitrary distribution was considered in  \cite{Yafaev3}. Then the role of the measure $d\Sigma (\z)$ is played also by a distribution which was called the  sigma-function of the Hankel operator $\bf H$.

\medskip
   
 {\bf 3.4.}  
  Although quite simple,  Theorem~\ref{HQ} is very useful for relating the discrete and continuous representations. Let us now discuss the case of bounded Hankel operators. Let
  $  D(z_{0},r)=\{|z -z_{0}|<r\}$  be the disc in the complex plane $\Bbb C$.   Recall that $dM(z)$ is called the  {\it Carleson measure}  on the unit disc $\Bbb D$ if  
   \begin{equation}
  \sup_{\mu  \in{\Bbb T}, r\in (0, r_{0})}  r^{-1}|M| (D(\mu,r)\cap{\Bbb D})<\infty  
  \label{eq:CA1}\end{equation}
  where   $r_{0} $ is some fixed small number.   
  Similarly, $d\Sigma(\z)$ is called the  Carleson measure  on the right half-plane ${\Bbb C}^+$ if  
    \begin{equation}
  \sup_{\lambda \in {\Bbb R}, R>0} R^{-1} |\Sigma| (D(i \lambda ,R)\cap{\Bbb C}^+)<\infty  .
 \label{eq:CA2}\end{equation}

    \begin{theorem}\label{Carle} 
    Let the measures  $dM(z)$  and $d \Sigma (\z)$ be linked by  the equation \e{eq:CM3}. Then
  $dM(z)$ is a Carleson measure on the unit disc $\Bbb D$ if and and only if $d \Sigma (\z)$   is a Carleson measure on the half-plane $\Bbb C^+$.
    \end{theorem}
    
   This assertion is checked in the Appendix by  straightforward but rather tedious calculations.  
    Theorem~\ref{Carle}  can also be indirectly deduced from general results on analytic functions (see, e.g., Section~E in Chapter~VIII of   \cite{Koosis}).
    
     \begin{theorem}\label{CD}  
   A Hankel operator $  H$   in the space $\ell^2 ({\Bbb Z}_{+})$
     is bounded  if and only if its matrix elements admit the representation \e{eq:CM1} where $dM(z)$ is some  Carleson measure on the unit disc $\Bbb D$.
    \end{theorem}
  
    This result is stated as Theorem~7.4 in Chapter~1 of the book \cite{Pe} and can be easily deduced from  Theorem~A2.12 of  \cite{Pe}
  where the representation of the symbol of a bounded Hankel operator as the Poisson balayage of some Carleson measure is given. For the proof of the latter result, we refer to Section~G in Chapter~X of the book \cite{Koosis}  or pages 271, 272 of the book \cite{Garn}.
   
   Putting together Theorems~\ref{HQ}, \ref{CD} and Lemma~\ref{Carle}, we get   the following result.
  
   \begin{theorem}\label{ex1}  
   A Hankel operator $\bf H$   in the space $L^2 ({\Bbb R}_{+})$
     is bounded if and only if its kernel  admits the representation \e{eq:CM2} with    some   Carleson measure $d\Sigma(\z)$    on the right half-plane ${\Bbb C}^+$.
    \end{theorem}
    
   It  also follows from Corollary~\ref{HQb} that if $dM(z)$ and $d\Sigma (\z)$ in the equation \e{eq:CM3} are  Carleson measures, then   the corresponding operators $H$ and $\bf H$ are unitarily equivalent.
     
    Of course  the representations \e{eq:CM1} and \e{eq:CM2} are highly non-unique. Let us give a simple 
    
    \begin{example}\label{uniq}  
    Let $d {\bf M}_{0}(z) $ be the Lebesgue measure on $\Bbb D$,  and let the measure $d N(z)$ be supported by the point $0$ so that $N(\{0\})=1$ and $N({\Bbb D}\setminus\{0\})=0$. Put 
    \[
    dM(z)= d {\bf M}_{0} (z) -\pi dN (z).
    \]
    Then
    \[
    \int_{\Bbb D} z^n dM(z)=   \int_{\Bbb D} z^n d {\bf M}_{0}(z)=\int_{0}^1 dr r^{n+1} \int_{0}^{2\pi}e^{in\theta} d\theta=0\q {\rm if}\q n\geq 1,
    \]
    and
    \[
    \int_{\Bbb D}   dM(z)=   \int_{\Bbb D}   d{\bf M}_{0}(z)- \pi=0.
    \]
     \end{example}
     
      \medskip
   
 {\bf 3.5.}
   In particular, the measures $dM(z)$ and $d\Sigma(\z)$ may be carried by the intervals $( -1,1)$ and $( 0,\infty)$, respectively.  Then  conditions  \e{eq:CA1} and  \e{eq:CA2} mean that
         \begin{equation}
   | M| ((1-\varepsilon,1))= O(\varepsilon)\q {\rm and} \q | M |((-1, -1+\varepsilon))= O(\varepsilon)\q {\rm as} \q \varepsilon\to 0
    \label{eq:pos}\end{equation}
    and
        \begin{equation}
  \sup_{  R>0} R^{-1} |\Sigma| (0,R)<\infty  .
 \label{eq:posD}\end{equation}

 For non-negative Hankel operators, the results of the previous subsection can be stated in a simpler and more definite form.  Recall that the Hankel form $h[f,f]$ was defined by relation \e{eq:B2h}. The following assertion is the   classical result of H.~Widom.
           
              \begin{theorem}\label{Widom}\cite[Theorem~3.1]{Widom}
              Suppose that $h[f,f]\geq 0$ for all $f\in {\cal D}_{+}$.        Then the following conditions are equivalent: 
                   \begin{enumerate}[\rm(i)]
\item
The operator $H$ is bounded.   
   \item
   The representation \e{eq:W} holds with a non-negative measure $dM(\nu)$ satisfying
   condition \e{eq:pos}.
  \item
  $a_{n}=O(n^{-1})$ as $n\to \infty$.
\end{enumerate}
 \end{theorem} 

Let us now state  a continuous analogue of Theorem~\ref{Widom}.

\begin{theorem}\label{WC} 
Let ${\bf a}\in C_0^{\infty} ({\Bbb R}_{+})'$.              Suppose that ${\bf h}[ {\bf f}, {\bf f}]\geq 0$ for all ${\bf f}\in C_0^{\infty} ({\Bbb R}_{+})$.        Then the following conditions are equivalent: 
                   \begin{enumerate}[\rm(i)]
\item
The operator $\bf H$ is bounded.   
   \item
   The representation \e{eq:S} holds with a non-negative measure $d\Sigma(\lambda)$   
   satisfying   condition \e{eq:posD}.
  \item
  ${\bf a} (t )=O(t^{-1})$ as $t\to \infty$ and $t\to 0$.
\end{enumerate}
 \end{theorem}
 
 \begin{pf}
 According to Theorem~5.1 in \cite{Yafaev3} the condition ${\bf h}[ {\bf f}, {\bf f}]\geq 0$ implies that, for some non-negative measure $d\Sigma(\lambda)$ on   ${\Bbb R} $,
  \begin{equation}
{\bf a}(t) = \int_{-\infty}^\infty e^{-\lambda t} d\Sigma(\lambda) 
 \label{eq:SS}\end{equation}
 where the integral converges for all $ t>0$. If $\Sigma({\Bbb R} \setminus {\Bbb R}_{+}) > 0$, then ${\bf a}(t) \geq c>0$ so that the operator $\bf H$ cannot be bounded. Thus representation \e{eq:SS} reduces to \e{eq:S}. 
 Let the measures $d\Sigma(\lambda) $ and $ dM(\nu) $ be linked by equation \e{eq:MS}.  By Lemma~\ref{Carle} (which is quite easy in this particular case) the conditions \e{eq:pos} and \e{eq:posD} are equivalent and hence, by Theorem~\ref{HQ}, the operator ${\bf H}$ is bounded if and only if  the operator $H$   with matrix elements \e{eq:W} is bounded. So the conditions (i) and (ii) are equivalent according to Theorem~\ref{Widom}.
 
It follows from  \e{eq:SS} that ${\bf a}\in C^{\infty} ({\Bbb R}_{+})$. So under assumption (iii) the operator ${\bf H}$ is bounded because the Carleman operator (it has integral kernel ${\bf a} (t )= t^{-1}$)  is bounded. Conversely, integrating in \e{eq:S} by parts we see that
  \[
{\bf a}(t) = t \int_{0}^\infty e^{-\lambda t} \Sigma(\lambda) d \lambda. 
 \]
 Therefore    condition \e{eq:posD} implies (iii).
\end{pf}

 Of course under the assumptions of Theorems~\ref{Widom} and \ref{WC} the measures $dM(\nu)$  and
  $d\Sigma(\lambda)$ are unique. 
     
 \medskip
   
 {\bf 3.6.} 
      By their definitions,  Carleson measures are carried by the open sets $\Bbb D$ or ${\Bbb C}^+$. Let us now consider the opposite case when  the measure $dM(z)$ is supported
   on the unit circle $\Bbb T$ and the corresponding measure $d\Sigma(\z)$ is supported
   on the line $\Re \z =0$. Then relations \e{eq:CM1}, \e{eq:CM2} and \e{eq:CM3} read, respectively, as
  \begin{equation}
a_{n } = \int_{ \Bbb T}\mu^n dM(\mu), \q n=0,1,\ldots, 
 \label{eq:CM1T}\end{equation}
   \begin{equation}
{\bf a}(t) = \int_{ \Bbb  R} e^{-i\lambda t} d\Sigma(i \lambda), \q  t>0,
 \label{eq:CM2T}\end{equation}
 and
   \begin{equation}
d M(\mu)= -2\alpha (\lambda-i \alpha)^{-2} d\Sigma (i \lambda), \q \mu=\frac{\lambda+i \alpha}{\lambda-i \alpha}.
 \label{eq:CM3T}\end{equation}
 
 In particular, if the measures $dM(\mu)$ on $\Bbb T$ and $d\Sigma(i \lambda)$ on $\Bbb R$ are absolutely continuous, that is,   
$dM(\mu)= \omega(\bar{\mu}) d {\bf m}_{0}(\mu) $ and $d\Sigma (i\lambda)=(2\pi)^{-1} \Omega (\lambda) d\lambda$, then \e{eq:CM1T}, \e{eq:CM2T} give the representations of the matrix elements $a_{n}$ of $H$ and the integral kernel ${\bf a}(t)$ of $\bf H$ in terms of their symbols $\omega (\mu)$ and $\Omega(\lambda)$ (see  formulas \e{eq:QF1}  and \e{eq:QF1c}). In this case \e{eq:CM3T} yields the standard relation \e{eq:B5h} 
 between these symbols. We recall that by the Nehari  theorem \cite{Nehari}, the operators $H$ or $\bf H$ are bounded if and only if the symbols $\omega$ or $\Omega$ can be chosen as bounded functions, but the construction above does not require this condition.

\section{Singular case}


 {\bf 4.1.}   
In this section, we consider   
 {\it signed } real measures $d\Sigma (\lambda)$ and $dM(\nu)$ satisfying the assumptions of Theorem~\ref{HQX}, but we do not assume that 
they satisfy the Carleson conditions \e{eq:pos} or \e{eq:posD}. We suppose that  these measures are absolutely continuous.
 Then \e{eq:W} and  \e{eq:S}  read  as
   \begin{equation}
                    {\bf a}(t) = \int_{0}^\infty e^{-\lambda t}  \sigma(\lambda) d\lambda
 \label{eq:W1}\end{equation}
 and
   \begin{equation}
   a_{n } = \int_{-1}^1 \nu^n \eta (\nu) d\nu , \q n=0,1,\ldots, 
 \label{eq:HD4}\end{equation}
 with some real functions $\sigma $ such that
  \[
  \int_{0}^\infty (\lambda+1)^{-2}  |\sigma(\lambda) |d\lambda <\infty
 \]
 and $\eta\in L^1  (-1,1)$. The relation \e{eq:MS} yields
  \begin{equation}
 \eta (\nu)=\sigma \big( \alpha\frac{ 1+\nu}{ 1-\nu}\big)  .
 \label{eq:MS1xld}\end{equation}
 
  It turns out that even this particular case leads to interesting examples. 
    Our plan is to use the results obtained in \cite{Yf1a, Ya} for integral Hankel operators $\bf H$, and then with the help of 
  Theorem~\ref{HQX} to state new results for  matrix Hankel operators $ H$. 
  
  Let us first illustrate our approach on  the Carleman operator and the Hilbert matrix already discussed in Subsection~1.3. If   $\sigma (\lambda)=1$, then according to \e{eq:W1} we have ${\bf a}(t)=t^{-1}$       which yields the Carleman operator $\bf H$.  It follows from  \e{eq:MS1xld}    that
$\eta (\nu)=1$ and hence, by  \e{eq:HD4}, the corresponding ``discrete" Hankel  operator $H ={\cal L}^* {\bf H}{\cal L} $ has matrix elements
\[
   a_n=(1+ (-1)^n) (n+1)^{-1}.
   \]

Let $\1_{X}$ be the characteristic function of a set $X\subset \Bbb R$.  Suppose that   $\eta (\nu)=\1_{(0,1)} (\nu)$.
Then according to  \e{eq:HD4} we have $a_n=  (n+1)^{-1}$
which yields the Hilbert matrix $H$. It follows from \e{eq:MS1xld} that
  $\sigma (\lambda)=\1_{(\alpha,\infty)} (\lambda)$ so that  the corresponding ``continuous" Hankel  operator ${\bf H} ={\cal L}  H {\cal L}^* $   has integral  hernel ${\bf a}(t)=t^{-1}e^{-\alpha t}$.

  Our goal is to study the case when the sigma-function $\eta (\nu) $ is singular at the points $\nu=\pm 1$ and the Hankel  operator $H$ with matrix elements \e{eq:HD4}
 is unbounded. We will consider two families of such Hankel  operators. The first family admits an explicit spectral analysis (Theorem~\ref{SpThHa}). Spectral information about the second family is more limited (Theorem~\ref{QC}).
 
 \medskip
 
 {\bf 4.2.}  
First we consider functions $\eta (\nu) $ with  arbitrary  logarithmic singularities at the points $\nu=  1$ and $\nu=- 1$.  To be   precise, we assume that
       \begin{equation}
\eta(\nu)=\sum_{l=0}^p
 \gamma_{l} \ln^l \big(\alpha\frac{1+\nu}{1-\nu}\big), \q p\geq 1,
 \label{eq:HD3}\end{equation}
 where $  \gamma_l$, $l=0,1,\ldots, p$, are any real numbers. Without loss of generality, we set $ \gamma_p=1$.   According to   \e{eq:HD4} the matrix elements of the corresponding Hankel operator \e{eq:B1h} are given by the formulas
     \begin{equation}
a_{n}=\sum_{l=0}^p
 \gamma_{l} \int_{-1}^{1 }\nu^{n}\ln^l \big(\alpha\frac{1+\nu}{1-\nu}\big)d\nu.
 \label{eq:HD3M}\end{equation}
 Note  (see  Proposition~\ref{SpThHx}, below) that
  \begin{equation}
a_{n}
 = ( 1  + (-1)^ {n+p} )  \frac{ \ln^p n }{n}
   \big(1+ O (  \frac{1 }{\ln n} )\big)
 \label{eq:DR7}\end{equation}
as $n\to\infty$.  It can be expected that such Hankel operator $H$ is unbounded because $\eta(\nu) d\nu$ is not the Carleson measure and $| a_{n}| n\to\infty$ as $n\to\infty$. Nevertheless, since  $a_{n} =\bar{a}_{n}$ and $\{a_{n}\}_{n\in {{\Bbb Z}_{+}}}\in \ell^{2}({\Bbb Z}_{+})$, the operator $H$ is well defined and symmetric on the dense set 
${\cal D}_{+} \subset \ell^{2}({\Bbb Z}_{+})$ of sequences with only a finite number of non-zero components.

Our study of the Hankel  operator $H$ with  matrix elements   \e{eq:HD3M} relies on  a combination of Theorem~\ref{HQX} with the results of \cite{Yf1}, \cite{Ya}
on integral Hankel operators $\bf H$ in the space $L^2 ({\Bbb R}_{+})$. In view of \e{eq:MS1xld} the sigma-function of ${\bf H}={\cal L} H {\cal L}^{*}$ equals
\[
\sigma(\lambda)=\sum_{l =0}^p
  \gamma_l \ln^l  \lambda, \q p\geq 1,
\]
and its integral kernel ${\bf a}(t)$   is given by formula \e{eq:W1}.  Let us define
a differential operator 
        \begin{equation}
L= v \sum_{l =0}^p
  \gamma_l D^l  v ,\q D    = i d /d \xi,
\label{eq:LOGz}\end{equation}
 of order $p$ in the space $L^2 ({\Bbb R}) $
where $v$ is the operator of multiplication by the universal function 
  \[
v (\xi )=\frac{\sqrt{\pi}} {\sqrt{\cosh (\pi \xi)}} .
\]
Let us introduce also a unitary transformation $F: L^{2}({\Bbb R}_{+})\to L^{2}({\Bbb R} )$ by the formula
\[
(F{\bf f}) (\xi)= (2\pi)^{-1/2}\frac{\Gamma (1/2 + i\xi)} {|\Gamma (1/2 + i\xi)|}\int_{0}^{\infty}t^{-1/2 -i\xi} {\bf f} (t) dt
\]
(so it is the Mellin transform, up to an insignificant phase factor). Note that $F{\bf f}\in{\cal S}$ if ${\bf f}\in{\bf D}_{+}
= {\cal L} {\cal D}_{+}$.

Putting together Theorem~\ref{HQX} and Theorem~3.2 of \cite{Yf1}, we can state the following result.
       
 \begin{proposition}\label{Qua}
Let matrix elements of the Hankel   operator $H$ be defined by formula \e{eq:HD3M} where the function $\eta(\nu)$is given by formula \e{eq:HD3}. Then for all $f\in {\cal D}_{+}$, the identity
\[
(Hf,f)= (LF {\cal L} f,F{\cal L} f).
\]
holds.
    \end{proposition}
    
It is shown in     Theorem~3.13 of \cite{Yf1} that the operator $L$ defined  on the set $C_{0}^{p}({\Bbb R}_{+})$ is essentially self-adjoint. The same arguments work  if $L$ is defined on the set $F  {\bf D}_{+}$.
 Therefore Proposition~\ref{Qua} allows us to obtain a similar result for the operator $H$.

 \begin{proposition}\label{Qua1}
 The Hankel   operator $H$ with matrix elements  \e{eq:HD3M} is essentially self-adjoint on the set  ${\cal D}_{+}$.
    \end{proposition}

Our goal 
is to obtain rather a complete information about the spectral structure of the closure of $H$    which will  also be  denoted $H$.

  \begin{theorem}\label{SpThHa}
Let matrix elements of the Hankel   operator $H$ be defined by formula \e{eq:HD3M} where $\gamma_{p}=1$ and $\gamma_{\ell} $ for $\ell \leq p-1$ are arbitrary real numbers.  
Then
\begin{enumerate}[\rm(i)]
\item
The spectrum of the operator $H$ is absolutely continuous except eigenvalues that  may accumulate to zero and  infinity only.   
   \item
The 
 absolutely continuous  spectrum of the operator $H$ covers $\Bbb R$ and is simple for $p$ odd. It coincides with $[0,\infty)$ and has multiplicity $2$  for $p$ even. 
 \item
 The point $0$ is not an eigenvalue of the operator $H$.
 If $p$ is odd, then the multiplicities of eigenvalues  of the operator $H$  are bounded by $(p-1)/2$. If $p$ is even, then the multiplicities of positive eigenvalues are bounded by $p/2-1$, and the multiplicities of negative eigenvalues are bounded by $p /2$. 
\end{enumerate}
    \end{theorem}
    
    Given Proposition~\ref{Qua}, this result  is a consequence (except the assertion about the point $0$ which is proven in  Theorem~4.7 of \cite{Yf1})  of the corresponding statement, Theorem~4.8 in \cite{Ya}, for the differential operator $L$.  
We emphasize that the method of \cite{Ya} yields a sufficiently explicit spectral analysis of the  operator $L$ and, in particular, information about its eigenfunctions of the continuous spectrum. In view of the unitary equivalence,  this yields the corresponding results for the Hankel  operator   $H$ (and ${\bf H}={\cal L} H{\cal L}^*$), but we will not dwell upon them.

 It remains to justify asymptotic relation \e{eq:DR7}.
 
  \begin{lemma}\label{log}
The asymptotic relation 
    \begin{equation}
  \int_{0}^{1} \nu^{n} \ln^{l } (1-\nu)  d\nu=   \frac{(-\ln n)^l }{n}
\big( 1   + O ( \frac{1 }{\ln n}) \big).
 \label{eq:DR5}\end{equation}
as $n\to\infty$  holds.
       \end{lemma}

  \begin{pf}
   Differentiating the formula
   \begin{equation}
  \int_{0}^{1} \nu^{n}  (1-\nu)^{\epsilon} d\nu=\frac{\Gamma(n+1) \Gamma(\epsilon+1)}{  \Gamma(\epsilon+n+2)}
 \label{eq:DR3}\end{equation}
 $l$ times in $\epsilon$ and then putting $\epsilon=0$, we see that
  \begin{equation}
  \int_{0}^{1} \nu^{n} \ln^{l } (1-\nu)  d\nu=   \frac{\partial^{l} }{\partial \epsilon^{l} }
  \Big( \frac{\Gamma(n+1)  \Gamma(\epsilon+1)}{  \Gamma(\epsilon+n+2)}\Big)\Big|_{\epsilon=0} .
 \label{eq:DR4}\end{equation}
 According to formula (1.18.4) in \cite{BE}  we have
   \begin{equation}
\frac{\Gamma(n+1)   }{  \Gamma(\epsilon+n+2)} = n^{-\epsilon-1} \big(1+O (n^{-1})\big),
 \label{eq:BE}\end{equation}
 and this asymptotic formula can be infinitely differentiated in $\epsilon$. In view of  \e{eq:DR4}, this yields  \e{eq:DR5}.
     \end{pf}

 Integrating separately over the intervals $(-1,0)$, $(0,1)$ and making  the change of variables $\nu\mapsto -\nu$
  in the  integral over   $(-1,0)$, we see that  
   \begin{align*}
 \int_{-1}^{1} \nu^n \ln^l \big(\alpha\frac{1+\nu}{1-\nu}\big) d\nu
 =  (-1)^l \int_0^{1} \nu^n \big( \ln(1-\nu) - \ln(1+\nu)- \ln \alpha\big)^l d\nu
\nonumber \\
 +  (-1)^{n} \int_0^{1} \nu^n \big( \ln(1-\nu) - \ln(1+\nu)+ \ln \alpha\big)^l d\nu.
 \end{align*}
 The asymptotic behavior   of the integrals on the right is determined by a neighborhood of the point $\nu=1$ where the function $\ln(1-\nu)$ is singular. The terms with $\ln(1+\nu)$ and $\ln\alpha$ do not give a contribution to the leading term
 of the asymptotics. Therefore putting together formulas  \e{eq:HD3M} (where $\gamma_p=1$) and \e{eq:DR5} we obtain  the following result.

  \begin{proposition}\label{SpThHx}
  Matrix elements \e{eq:HD3M} of the   Hankel operator  $H$ satisfy the asymptotic relation \e{eq:DR7}.
    \end{proposition}
    
    We finally note that for $p=1$, the differential operator  \e{eq:LOGz} reduces by an explicit unitary transformation to 
    the  operator $  i d /d \xi $. Therefore the same is true for  the corresponding Hankel operators $\bf H$ and $H$, see \cite{Yf1} for details.
    
\medskip
 
     {\bf 4.3.} 
  Now we consider even more singular  compared to \e{eq:HD3} case  when  the function $\eta (\nu) $ has power singularities at the points $\nu=  1$ or  $\nu= - 1$. Let
       \begin{equation}
\eta(\nu)=  \Big(\frac{1+\nu}{1-\nu}\Big)^q, \q q\in (-1, 1), \q q\neq 0,
 \label{eq:QC}\end{equation}
 and let the sequence $a_{n}$ be defined by formula \e{eq:HD4}. Then according to Theorem~\ref{HQXb} the non-negative Hankel operator $H$ with such matrix elements is correctly defined  via its quadratic form.
  Our goal here is to describe its spectral structure.

  \begin{theorem}\label{QC}
The spectrum of the Hankel operator $H=H(q)$ with the matrix elements
 \begin{equation}
 a_{n } = a_{n }(q) = \int_{-1}^{1} \nu^n \Big(\frac{1+\nu}{1-\nu}\Big)^q d\nu ,  \q q\in (-1, 1), \q q\neq 0,
 \label{eq:QC2}\end{equation}
 is absolutely continuous, coincides with the half-axis $[0,\infty)$ and has constant multiplicity. 
    \end{theorem}
    
  In view of  Theorem~\ref{HQXb}     this result can be deduced   from 
 the corresponding assertion for the Hankel operator ${\bf H}= {\bf H}(q)$ in the space $L^2 ({\Bbb R}_{+})$. Indeed, it follows from formula 
  \e{eq:MS1xld} that the sigma-function of the operator ${\bf H}={\cal L} H{\cal L}^*$ is
  $
  \sigma(\lambda)=2^q \lambda^q,
  $
  and hence according to relation   \e{eq:W1} its integral kernel equals
    \begin{equation}
  {\bf a}(t)=  {\bf a}_q(t)= 2^q \int_{0}^\infty e^{-t\lambda} \lambda^q d\lambda= 2^q \Gamma(q+1) t^{-q-1}.
  \label{eq:dila}\end{equation}
  Hankel operators $\bf H$ with such kernels were studied in \cite{Yf1a}. It was shown in  Theorem~1.2 that for all $q >-1$, $q\neq 0$, the spectra of the operators $\bf H$ are absolutely continuous, coincide with the half-axis $[0,\infty)$ and have constant multiplicity. In particular, for $q <1$, this result yields Theorem~\ref{QC}. 
  
  \begin{remark}\label{rmk1}
\begin{enumerate}[1.]
\item
The proof of Theorem~1.2 in \cite{Yf1a} relies only on the invariance of the Hankel operator with kernel \e{eq:dila} with respect to the group of dilations. So,  spectral information about the Hankel operators with matrix elements \e{eq:QC2} is more limited than under the assumptions of Theorem~\ref{SpThHa}. Even the spectral multiplicity of $\bf H$ is unknown. 
We recall that according to the fundamental results of \cite{MPT} the spectral multiplicity of a positive bounded Hankel operator does not exceed $2$, but, strictly speaking, this result is not applicable because the operator $H$ in Theorem~\ref{QC} is unbounded. Actually, we expect that the spectrum of $H$ is simple since the kernel \e{eq:dila} has only one singular point. This point is $t=0$ if $q>0$ and $t= \infty$ if $q <0$.

\item
It is by no means obvious how to prove Theorem~\ref{QC} directly in   $\ell^2 ({\Bbb Z}_{+})$ 
because the realization in this space of the group of dilations  in the space $L^2 ({\Bbb R}_{+})$ is  not transparent.

\item
 It follows from definition \e{eq:QC2} that
    \begin{equation}
 a_{n }(q) =(-1)^n  a_{n }(- q)
 \label{eq:U-}\end{equation}
 and hence  $ H(-q)= V^* H(q) V$ where the unitary operator $V$ is defined on sequences $f=(f_{0}, f_{1}, \ldots)$ by the formula $(V f)_{n} = (-1)^n f_{n} $. Thus the Hankel operators ${\bf H}(-q)$ and ${\bf H} (q)$ with kernels \e{eq:dila}  in the space $L^2 ({\Bbb R}_{+})$ are also unitarily equivalent. This fact does not look obvious in the continuous representation.

\item
For $q\geq 1$, equality \e{eq:QC2} makes no sense. In this case there is   no reasonable interpretation of   the Hankel operator $H$ with sigma-function \e{eq:QC}  in the  space $\ell^2 ({\Bbb Z}_{+})$ although the Hankel operator $\bf H$ with integral kernel \e{eq:dila} is well defined in the  space $L^2 ({\Bbb R}_{+})$. 
\end{enumerate}
\end{remark}

   Finally, we find the asymptotics of the matrix elements
 \e{eq:QC2}    as $n\to\infty$. In view of formula \e{eq:U-} we may suppose that $q\in (0,1)$.
 Then  the asymptotics  of  the integral  \e{eq:QC2} is determined by a neighborhood of the point $\nu=1$.    So,  we write formula  \e{eq:QC2} as
  \begin{multline*}
a_{n}=2^q \int_0^{1} \nu^n (1-\nu)^{-q} d\nu 
\\
+ \int_0^{1} \nu^n (1-\nu)^{-q}\big((1+\nu)^q -2^{q} \big) d\nu
+ (-1)^{n} \int_0^{1} \nu^n \Big(\frac{1-\nu}{1+ \nu}\Big)^q d\nu  .
\end{multline*}
 The first integral on the right  coincides with expression \e{eq:DR3}  for $\epsilon=-q$, and its asymptotics as $n\to\infty$ is given  by formula \e{eq:BE}.
 The second and third integrals are $O(n^{-2+q})$ and $O(n^{-1-q})$, respectively.
 This yields     the following result.

     \begin{proposition}\label{QC1}
   The sequence   \e{eq:QC2} satisfies the asymptotic relation 
     \[
a_{n}
 =  (\sgn q)^n 2^{|q|} \Gamma(1-|q|) n^{|q|-1}  \big(1+O (n^{-\varepsilon})\big), \q n\to \infty,
 \]
 where $\varepsilon= \min\{1, 2 |q|\}$.
    \end{proposition}
    
 As could be expected, $a_{n}\to 0$ as $n\to\infty$ but much slower than sequence  \e{eq:DR7}.

\appendix
\section{Proof of Theorem~\ref{Carle}}

{\bf A.1.}
Suppose, for definiteness,  that  the equations \e{eq:CM3} are satisfied with $\alpha=1$. Then 
\[
\zeta=\frac{1+z}  { 1-z}=:\omega(z)
  \q{\rm and}\q  z=\omega^{-1}(\zeta) =\frac{\zeta -1}  {\zeta +1}.
\]
A standard calculation shows that if  $r< |z_{0} -1|$, then
     \begin{equation}
 \omega(D(z_{0},r)) =D(\z_{0}, R)  
 \label{eq:FL}\end{equation}
where
   \begin{equation}
 \z_{0}=\frac{(1+z_{0}) (1-\bar{z}_{0})+r^2}{|1-z_{0}|^2 -r^2}, \q R =  \frac{2r}{|1-z_{0}|^2 -r^2}.
 \label{eq:FL1}\end{equation}
Conversely, supposing that $R <|\z_{0}+1|$, we see that relation \e{eq:FL} holds true with
   \begin{equation}
 z_{0}=\frac{(\z_{0}- 1) ( \bar{\z}_{0}+1) -R^2}{|\z_{0}+1|^2 -R^2}, \q r=  \frac{2 R}{|\z_{0}+1|^2 -R^2}.
 \label{eq:FL2}\end{equation}
 
 In particular,  if $z_{0}=e^{i\theta}$, $\theta\in [0,2\pi)$,   then $  |z_{0} -1|=2 \sin  (\theta /2) $ and \e{eq:FL1} reads as
    \begin{equation}
 \z_{0}=\frac{2i \sin\theta +r^2}{4\sin^{2} (\theta/2)  -r^2}, \q  R=  \frac{2r}{4\sin^{2} (\theta /2)  -r^2}.
 \label{eq:FL1a}\end{equation} 
 Similarly, if $\z_{0}= i\lambda$, $\lambda\in {\Bbb R}$, then $  |\z_{0} +1|=\sqrt{\lambda^2 +1} $ and \e{eq:FL2} reads as
 \begin{equation}
 z_{0}=\frac{\lambda^2 +2i\lambda -1 -R^2}{\lambda^2 +1- R^2}, \q r=  \frac{2 R}{\lambda^2 +1- R^2}.
 \label{eq:FL2a}\end{equation}

\medskip

{\bf A.2.}
By the proof of Theorem~\ref{Carle},
  we may suppose that the measures  $dM(z)$  and $d \Sigma (\z)$ are non-negative. It  is also convenient to extend these measures onto the whole complex plane ${\Bbb C}$ setting $M({\Bbb C}\setminus {\Bbb D})=0$ and $\Sigma({\Bbb C}\setminus {\Bbb C}_{+})=0$. Below $C$ (sometimes with indices) are different positive constants whose precise values are of no importance.
  
   Let us define a measure $d\wt{M}(\zeta)$  on ${\Bbb C}$ by the relation $\wt{M}(  Y)= M( \omega^{-1} (Y))$ 
 for an arbitrary set $Y\subset {\Bbb C}$. In particular, according to   \e{eq:FL} we have
    \begin{equation}
M(D(z_{0},r)) =\wt{M} (D(\z_{0}, R)  ).
 \label{eq:FLX}\end{equation}
  It follows from \e{eq:CM3}    that  
  \begin{equation}
 d\wt{M}(\zeta)=2 |\z+1|^{-2}d \Sigma (\z).
  \label{eq:A2}\end{equation}
  
  It is convenient to state the conditions \e{eq:CA1} and \e{eq:CA2} in an equivalent way.
  
   \begin{lemma}\label{AA1}
     \begin{itemize}
         \item
    Condition \e{eq:CA1} is satisfied if 
     \begin{equation}
  \sup_{\theta\in [0,2\pi)}  \sup_{  r \leq 2\gamma_{0} \sin (\theta/2) }  r^{-1}M (D(e^{i\theta},r) )<\infty   
  \label{eq:CA1X}\end{equation}
for some $\gamma_{0} \in (0,1)$ and, for some    $r_{0}>0$,
   \begin{equation}
  \sup_{ r\in (0, r_{0})}  r^{-1}M (D(1,r) )<\infty   .  
    \label{eq:CA1Y}\end{equation}
      \item
    Condition \e{eq:CA2} is satisfied if 
     \begin{equation}
  \sup_{ \lambda\in {\Bbb R}}    \sup_{ R \leq \d_{0} \sqrt{\lambda^2+1} } R^{-1} \Sigma (D(i \lambda ,R) )<\infty  
 \label{eq:CA2X}\end{equation}
for some $\d_{0}\in (0,1)$     and,  for some $R_{0} >0$,
    \begin{equation}
  \sup_{R \geq R_{0}} R^{-1} \Sigma (D(0 ,R) )<\infty  .
 \label{eq:CA2Y}\end{equation}
     \end{itemize}
    \end{lemma}
    
      \begin{pf}
   Indeed, if $ r \geq 2\gamma_{0} \sin (\theta/2) $, then $D(e^{i\theta},r) \subset D(1, \gamma r)$ for $\gamma =1+\gamma_0^{-1}$. Therefore \e{eq:CA1} for such $\theta$ and $r$ follows from
    \e{eq:CA1Y}.
    
   Similarly, if  $R\geq \d_{0}\sqrt{\lambda^2+1}$, then  $R\geq \d_{0}$ and $D(i\lambda, R) \subset D(0, \d R)$  for $\d =1+\d_{0}^{-1}$. Therefore \e{eq:CA2} for such $\lambda$ and $R$ follows from    \e{eq:CA2Y}.
   \end{pf}

First, we consider the conditions \e{eq:CA1X} and \e{eq:CA2X}.

  \begin{lemma}\label{A3}
   Condition  \e{eq:CA2} on $d\Sigma$ implies condition   \e{eq:CA1X} on $dM$.
       \end{lemma}
       
           \begin{pf}
    Let $r \leq 2\gamma_0 \sin (\theta/2)$, $\lambda=\cot (\theta/2)$,  and let $\z_{0}$ and $R$ be  given by formulas \e{eq:FL1a}.     
    Since $i\lambda\in  D(\z_{0}, R)$, we have
           $         D(\z_{0}, R)\subset  D(i\lambda  , 2 R)$,
      and hence      according to \e{eq:CA2}  
               \begin{equation}
         \Sigma(  D(\z_{0}, R ))\leq C R.
         \label{eq:AB2}\end{equation} 
         Put $z_{0}=e^{i\theta}$, then $ |1-z_{0}|=2\sin(\theta/2)$.
       If $\z\in  D(\z_{0}, R)$, then   $z= \omega^{-1} (\z)\in  D(z_{0}, r)$ and hence $|z_{0}-z|\leq r\leq 2\gamma_{0}\sin(\theta/2)$. So, we have
         \[
        2|1+\z|^{-1}=  |1-z|\leq  |1-z_{0}|+ |z_{0}-z| \leq 2 (1+\gamma_0) \sin(\theta/2).
\]
 In view of \e{eq:A2} it now follows from \e{eq:AB2} that
   \begin{equation}
         M(  D(z_{0}, r))=  \wt{M}(  D(\z_{0}, R))\leq  2  (1+\gamma_0)^2 \sin^2(\theta /2)  \Sigma(  D(\z_{0}, R))\leq C R\sin^2(\theta /2) .
         \label{eq:AB4}\end{equation} 
         Since, according to  the second formula  \e{eq:FL1a},
         \[
          R\sin^2(\theta /2) = 2^{-1} r \Big( 1-\frac{r^2}{4\sin^2(\theta/2)}\Big)^{-1}\leq  2^{-1}  \big( 1-\gamma_{0}^2\big)^{-1} r,
          \]
          estimate \e{eq:AB4} yields  \e{eq:CA1X}.
                        \end{pf}
                        
                        Next, we prove the converse assertion.
              
              \begin{lemma}\label{A4}
   Condition  \e{eq:CA1} on $d M$ implies condition   \e{eq:CA2X} on $d\Sigma$.
       \end{lemma}
       
           \begin{pf}
           Consider the disc $D(i\lambda,R)$ where   $R\leq \d_0\sqrt{\lambda^2+1} $ with $\d_0<1$. In view of      relations \e{eq:FLX} and \e{eq:A2}, we     have
            \begin{equation}
           \Sigma (D(i\lambda,R))  \leq C (\lambda^2+1)        \wt{M}(D(i\lambda,R)) =C (\lambda^2+1)       M(D(z_{0},r))                 \label{eq:D}\end{equation} 
                where $z_{0} $ and $r$ are given by  formulas \e{eq:FL2a}.
                Since $e^{i\theta}= z_{0} |z_{0} |^{-1}\in D(z_{0},r)$ we see that  $D(z_{0},r)\subset D( e^{i\theta} ,2r)$. Thus, it follows from condition \e{eq:CA1} that the right-hand side of \e{eq:D} is bounded by 
                \[
     C (\lambda^2+1)r =2     C R\Big(1- \frac{R^2}{\lambda^2+1}\Big)^{-1}\leq 2 \big(1- \d^2_{0}\big)^{-1}C  R.
                \]
 Therefore \e{eq:D} implies condition   \e{eq:CA2X}.
              \end{pf}
              
              \medskip

{\bf A.3}.
  It remains to compare the conditions \e{eq:CA1Y} and \e{eq:CA2Y}.
  
    \begin{lemma}\label{A2}
   Conditions  \e{eq:CA1Y} and \e{eq:CA2Y} are equivalent.
       \end{lemma}
       
           \begin{pf}
        Let \e{eq:CA2Y} be satisfied.
      Observe that
  \begin{equation}
 \omega(D(1,r)) ={\Bbb C}\setminus D(-1,2 r^{-1})
  \label{eq:A3}\end{equation}
and set 
$
\phi (R)= \Sigma (  D(-1 , R)    )  .
$
  In view of \e{eq:A2} we have
     \[
\wt {M}(  D(-1 ,R')  \setminus D(-1 ,R) )=\int_{R}^{R'} d \wt {M}(   D(-1 ,\rho)    )= 
2 \int_{R}^{R'} \rho^{-2} d \phi (\rho).
 \]
Integrating here by parts  and then   passing  to the limit $R'\to\infty$, we see that 
    \[
2^{-1}\wt {M}( {\Bbb C} \setminus D(-1 ,R))= - R^{-2} \phi (R)
+ 2 \int_{R}^{\infty}\rho^{-3}  \phi (\rho) d\rho.
 \] 
 Under assumption \e{eq:CA2Y}  the right- and hence the left-hand sides here are $O(R^{-1})$ as $R\to\infty$. Thus,
 according to  \e{eq:FLX} and \e{eq:A3}, we have
 \[
M (D(1,r))=\wt {M}( {\Bbb C} \setminus D(-1 ,2 r^{-1}))= O(r), \q r\to 0.
\]

 Conversely, let \e{eq:CA1Y} be satisfied.
   Again in view of \e{eq:A2} we have
   \[
2 \Sigma(   D(-1 ,R)  )=
 2 \int_{0}^R   d  \Sigma(   D(-1 , \rho) ) =   \int_{0}^R  \rho^{2} d \wt {M} (   D(-1 , \rho)  ) .
 \]
 Since, by \e{eq:A3},
  \[
  \wt {M} (  {\Bbb C} \setminus D(-1 , R))  =M (   D(1 ,2  R^{-1}))  =: \psi (R),
 \]
 it follows that
    \[
 2\Sigma(    D(-1 ,R)  )=-\int_{0}^R  \rho^{2} d \psi (\rho)= -R^{2} \psi (R)
+ 2 \int_{0}^R \rho   \psi (\rho) d\rho.
 \]
 Under assumption \e{eq:CA1Y}, $\psi (R)=O (R^{-1})$  as $R\to\infty$. So, the right- and hence the left-hand sides here are $O(R)$ as $R\to\infty$  which proves \e{eq:CA2Y}.
     \end{pf}

 Now it easy to conclude the proof of Theorem~\ref{Carle}. Let condition \e{eq:CA2} on $d\Sigma$ be true. Then, by
 Lemma~\ref{A3}, $dM$ satisfies \e{eq:CA1X} and, by Lemma~\ref{A2}, it satisfies \e{eq:CA1Y} . So, Lemma~\ref{AA1} implies that condition  \e{eq:CA1} is fulfilled. Similarly, let condition \e{eq:CA1} on $d M$ be true. Then, by
 Lemma~\ref{A4}, $d\Sigma$ satisfies \e{eq:CA2X} and, by Lemma~\ref{A2}, it satisfies \e{eq:CA2Y}. So, Lemma~\ref{AA1} implies that condition  \e{eq:CA2} is fulfilled. $\qed$



\begin{thebibliography}{00}



\bibitem{Bo} A.~B\"ottcher, B.~Silbermann,
\emph{Analysis of Toeplitz operators},  Springer-Verlag,  2006.

 



\bibitem{BE} A.~Erd\'elyi,  W.~Magnus, F.~Oberhettinger, F.~G.~Tricomi,
\emph{Higher transcendental functions}, Vol. 1, 2, 
McGraw-Hill, New York-Toronto-London, 1953.


\bibitem{Garn} J.~B.~Garnett,
\emph{Bounded analytic functions,  }  
Academic Press, New York,  1981.

\bibitem{GUEVI+} I.~M.~Gel'fand,  G.~E.~Shilov,  
\emph{Generalized functions.  } Vol.~1, 
Academic Press, New York and London,  1964.   

\bibitem{GGK}
 I.~Gohberg, S.~Goldberg, M.~Kaashoek,
\emph{Classes of linear operators},   vol. 1, Birkh\"auser,  1990. 


\bibitem {Hof} K. Hoffman, {\it Banach spaces of analytic functions}, Prentice-Hall, Inc., Englewood Cliffs, New York, 1962.

\bibitem{Koosis}
P.~Koosis,
\emph{Introduction to $H_{p}$ spaces}, 
LMS, Lecture Note Series {\bf 40},  Cambridge Univ. Press, 1980. 

\bibitem{Ma}
  W.~Magnus, 
\emph{On the spectrum of Hilbert's matrix}, Amer. J. Math. 
{\bf 72} (1950), 405-412.

\bibitem{MPT}
 A.~V.~Megretskii, V.~V.~Peller,  S.~R.~Treil,
\emph{The inverse spectral problem for self-adjoint Hankel operators}, Acta Math.   {\bf 174} (1995), 241-309. 

 
\bibitem{Nehari}
Z.~Nehari,  
\emph{On bounded bilinear forms}, Ann. Math.   {\bf 65} (1957), 153-162. 
 

 


 \bibitem{NK} N.~K.~Nikolski, {\em Operators, functions, and systems: an easy reading}, vol. I: Hardy, Hankel, and Toeplitz, Math. Surveys and Monographs vol.~92,  Amer. Math. Soc.,   Providence,
  Rhode Island, 2002.
  
 
  
\bibitem{Pe}
V.~V.~Peller,
\emph{Hankel operators and their applications}, 
Springer Verlag,  2002. 


   \bibitem{Ro}
  M.~Rosenblum, 
\emph{On the   Hilbert matrix}, I, II, Proc. Amer.   Math. Soc.
{\bf 9} (1958), 137-140, 581-585.

   \bibitem{RoRo}
  M.~Rosenblum,   J.~Rovnyak, 
\emph{Hardy classes and operator theory}, Oxford Univ. Press, 
1985.



\bibitem{Widom}
H.~Widom,
\emph{Hankel matrices}, 
Trans. Amer. Math. Soc. {\bf 121}  (1966), 1-35.


  

\bibitem{Yf1}
D.~R.~Yafaev, \emph{Diagonalizations of two classes of unbounded Hankel operators},   Bulletin Math. Sciences  {\bf 4} (2014), 175-198.

\bibitem{Y}
D.~R.~Yafaev, \emph{Criteria for Hankel operators to be sign-definite}, Analysis $\&$ PDE  {\bf 8} (2015), 183-221.



\bibitem{Yf1a}
D.~R.~Yafaev, \emph{Quasi-Carleman operators and their spectral properties},  Integral Eq. Op. Th.   {\bf 81} (2015), 499-534.

\bibitem{Yfr}
D.~R.~Yafaev, \emph{On finite rank Hankel operators},  J. Funct. Anal.   {\bf 268} (2015), 1808-1839.

\bibitem{Yafaev3}
D.~R.~Yafaev,
\emph{Quasi-diagonalization of Hankel operators,}
preprint, arXiv:1403.3941(2014); to appear in  J. d'Analyse Math\'ematique.

\bibitem{Ya}
D.~R.~Yafaev, \emph{ Spectral and scattering theory for differential  and Hankel operator}, arXiv: 1511.04683 (2015).

\bibitem{Yunb}
D.~R.~Yafaev, \emph{ Unbounded Hankel operator and moment problems},   Integral Eq. Oper. Theory, {\bf 85}  (2016), 289-300.

\bibitem{YaT}
 D.~R.~Yafaev, \emph{ On semibounded Toeplitz operators}, arXiv: 1603.06229 (2016), to appear in J. Oper. Theory.

\bibitem{YaWH}
D.~R.~Yafaev, \emph{ On semibounded Wiener-Hopf operators}, arXiv: 1606.01361 (2016).

\end{thebibliography}
 \end{document}